\RequirePackage{etoolbox}
\csdef{input@path}{{style/}{graphics/}}
\documentclass[aop,MSNbibl,dvips]{arximspdf}
\usepackage{graphicx}
\doi{10.1214/13-AOP884} 
\volume{43}
\issue{2}
\pubyear{2015}
\firstpage{605}
\lastpage{638}

\makeatletter
\newcommand{\rrVert}{\Vert}
\newcommand{\rrvert}{\vert}
\newcommand{\llVert}{\Vert}
\newcommand{\llvert}{\vert}
\newtheorem{lemma}{Lemma}[section]
\newtheorem{proposition}[lemma]{Proposition}
\newtheorem{theorem}[lemma]{Theorem}
\newtheorem{corollary}[lemma]{Corollary}
\newproclaim{remark}[lemma]{Remark}
\newproclaim{definition}[lemma]{Definition}
\newproclaim{example}[lemma]{Example}
\makeatother

\begin{document}
\begin{frontmatter}

\title{On fractional smoothness and $L_p$-approximation
on~the Gaussian space}
\pdftitle{On fractional smoothness and Lp-approximation
on~the Gaussian space}
\runtitle{Fractional smoothness}

\begin{aug}
\author[A]{\fnms{Stefan} \snm{Geiss}\corref{}\thanksref{T1}\ead[label=e1]{stefan.geiss@uibk.ac.at}}
\and
\author[B]{\fnms{Anni} \snm{Toivola}\thanksref{T2}\ead[label=e2]{anni.toivola@jyu.fi}}
\runauthor{S. Geiss and A. Toivola}
\affiliation{University of Jyv\"askyl\"a and University of Innsbruck,
and University of Jyv\"askyl\"a}
\address[A]{University of Jyv\"askyl\"a\\
P.O. Box 35\\
FIN-40014\\
Finland\\
and\\
University of Innsbruck\\
Technikerstra\ss e 13/7\\
A-6020 Innsbruck\\
Austria\\
\printead{e1}} 
\address[B]{University of Jyv\"askyl\"a \\
P.O. Box 35\\
FIN-40014\\
Finland\\
\printead{e2}}
\end{aug}
\thankstext{T1}{Supported by the Project 133914 Stochastic
and Harmonic Analysis, Interactions and Applications of the
Academy of Finland.}
\thankstext{T2}{Supported by the Emil Aaltonen Foundation and
the Magnus Ehrnrooth Foundation (Grant no.~MA2012n26).
Some part of the work was done as the author
was affiliated to the Department of Mathematics,
University of Innsbruck, Austria.}

\received{\smonth{10} \syear{2012}}
\revised{\smonth{8} \syear{2013}}

%
\begin{abstract}
We consider Gaussian Besov spaces obtained by real interpolation and
Riemann--Liouville operators of fractional integration on the Gaussian space
and relate the fractional smoothness of a functional to the regularity
of its heat extension. The results are applied to study an approximation
problem in $L_p$ for $2\le p < \infty$ for stochastic integrals with respect
to the $d$-dimensional (geometric) Brownian motion.
\end{abstract}

%

\end{frontmatter}

\section{Introduction}
This paper is devoted to Besov spaces defined on a Gaussian space, associated
Riemann--Liouville operators of fractional integration, and
approximation theory.
As Gaussian space, we consider $L_p(\mathbb{R}^d,\gamma_d)$ with
$2\le p < \infty$ and
$d\gamma_d = e^{-|x|^2/2} \,dx/(2\pi)^{d/2}$ being the
$d$-dimensional standard Gaussian measure.
The (Gaussian) Besov spaces are obtained by the real interpolation
method and the approximation problem
concerns an approximation of stochastic integrals in $L_p$. Some of the
results are extensions of
corresponding statements proved mainly in $L_2$; see \cite{ZhangR,Gobet-Temam,Gei15,Gei18,Gei17,Gei20,Hujo3,Hujo2007,Gei19,Seppala1}.
However, the $L_2$-theory and the $L_p$-theory
for $2<p<\infty$ on a Gaussian space may differ significantly. For
example, the Meyer inequalities
can be proved in $L_2$ using orthogonality by standard ideas, but they
are considerably more
involved in the $L_p$-case when $1<p\neq2<\infty$ (see \cite{Nualart2}, Proposition 1.5.3, and \cite{Pisier1988}). Another example
is the phenomenon that, for instance, for $2<p<\infty$
and $f\in L_p(\mathbb{R},\gamma_1)$, the orthogonal Hermite expansion
does not
necessarily converge in $L_p(\mathbb{R},\gamma_1)$ (see \cite{Pollard1948}).

Regarding the multi-step $L_p$-approximation problem on the Gaussian
space for $2<p<\infty$
we study in this paper, we cannot exploit chaos expansion techniques
like in
\cite{Geiss-Geiss-Laukkarinen2012} nor can we reduce the problem by
orthogonality to a question
about a one-step approximation as in the $L_2$-setting \cite{Gei15,Gei18} and
BMO-setting \cite{Gei17}.
The difference between the $L_2$- and the $L_p$-context for $2<p<\infty
$ is also visible
by the fact that we have to describe the optimal $L_p$-approximation in
Theorem~\ref{thmLp-adaptedtimenets} below by a Riemann--Liouville
operator instead of the real interpolation spaces.

To explain the purpose of this paper in more detail, let us introduce
some notation.
We let $W=(W_t)_{t\in[0,1]}$ be a standard $d$-dimensional Brownian motion
starting in zero defined on $(\Omega,\mathcal{F},\mathbb
{P},(\mathcal{F}_t)_{t\in[0,1]})$,
where $(\Omega,\mathcal{F},\mathbb{P})$ is complete and $(\mathcal
{F}_t)_{t\in[0,1]}$
is the augmentation of the natural filtration and where we can assume that
$\mathcal{F}= \mathcal{F}_1$. As processes driving the stochastic
integrals, we use
the Brownian motion and the coordinate-wise geometric Brownian motion,
that is,
\[
Y_t:=\bigl(W_t^{(1)},\ldots,W_t^{(d)}
\bigr)^\top\quad\mbox{and}\quad E:=\mathbb{R}^d
\]
or
\[
Y_t:=\bigl(e^{W_t^{(1)}-(t/2)},\ldots,e^{W_t^{(d)}-(t/2)}
\bigr)^\top\quad\mbox{and}\quad E:=(0,\infty)^d.
\]
Then we have
\[
dY_t = \sigma(Y_t) \,dW_t,
\]
where $Y$ is considered as a column vector and the $d\times d$-matrix
$\sigma(y)$
is given by $\sigma(y) = I_d$ or $(\sigma_{ij}(y))_{i,j=1}^d =
(\delta_{i,j} y_i)_{i,j=1}^d$, respectively,\vspace*{2pt}
where $\delta_{i,j}=1$ if $i=j$ and $\delta_{i,j}=0$ otherwise.
The parabolic differential operator associated to the diffusion~$Y$~is
\[
\mathcal{A}:= \frac{\partial}{\partial t} + \frac{1}{2} \sum
_{k=1}^d \sigma_{kk}^2
\frac{\partial^2}{\partial y_k^2}.
\]
Given a Borel-function $g\dvtx E\to\mathbb{R}$ with $g(Y_1)\in L_2$, we let
%
\begin{equation}
\label{eqndefinitionG} G(t,y):= \mathbb{E}\bigl(g(Y_1)| Y_t=y
\bigr)
\end{equation}
and notice that $G(1,y)=g(y)$. Integrability properties of $G$ and its
derivatives are
given in Lemma~\ref{lemmahypercontraction} below and are used
implicitly in this paper.
The function $G$ solves the backward parabolic PDE
\[
\mathcal{A} G = 0 \qquad\mbox{on }[0,1)\times E.
\]
For $0 \leq s < t < 1$, It\^o's formula implies that
%
\begin{equation}
\label{eqnItoG} G(t, Y_t) - G(s, Y_s) = \int
_s^t \nabla G (u, Y_u)
\sigma(Y_u) \,dW_u\qquad \mbox{a.s.},
\end{equation}
where $\nabla G(t,x)$ is considered as a row vector.
Furthermore,
%
\begin{equation}
\label{eqnrepresentationgY1} g(Y_1) = \mathbb{E}g(Y_1) + \int
_0^1 \nabla G(u, Y_u)
\sigma(Y_u) \,dW_u\qquad\mbox{a.s.}
\end{equation}
by $t\uparrow1$, where the convergence takes place in $L_2$ [or later in
$L_p$ if $g(Y_1)\in L_p$ with $2\le p < \infty$].
One purpose of this paper is to investigate Riemann approximations
of the stochastic integral in (\ref{eqnrepresentationgY1}) by
the following quantities.

%
\begin{definition}\label{deferrorprocesses}
(i)~Let ${\mathcal T}^{\mathrm{rand}}$ be the set of all sequences
of stopping times
$\tau=(\tau_i)_{i=0}^n$ with
$0=\tau_0 \le\tau_1 \le\cdots\le\tau_{n-1} <\tau_n = 1$ where
$n=1,2,\ldots,$ such that
$\tau_i$ is $\mathcal{F}_{\tau_{i-1}}$-measurable for $i=1,\ldots,n-1$,
that is,
\[
\{ \tau_i \in B \} \cap\{\tau_{i-1} \le t \} \in
\mathcal{F}_t \qquad\mbox{for } t\in[0,1] \mbox{ and } B\in
\mathcal{B}\bigl([0,1]\bigr).
\]

(ii)~Given a time-net $\tau=(\tau_i)_{i=0}^n\in{\mathcal
T}^{\mathrm{rand}}$, $0\le t\le
1$ and $g(Y_1)\in L_2$, we let
\begin{eqnarray*}
C_t\bigl(g(Y_1),\tau\bigr) &:= & \int
_0^t \nabla G(s,Y_s)
\,dY_s
- \sum_{i=1}^n \nabla G(
\tau_{i-1},Y_{\tau_{i-1}}) (Y_{\tau_i
\wedge t} - Y_{\tau_{i-1} \wedge t}),
\\
C_t\bigl(g(Y_1),\tau,v\bigr) &:= & \int
_0^t \nabla G(s,Y_s)
\,dY_s - \sum_{i=1}^n
v_{\tau_{i-1}} (Y_{\tau_i \wedge t} - Y_{\tau_{i-1}
\wedge t}),
\end{eqnarray*}
where $v = (v_{\tau_{i-1}})_{i=1}^n$ is a sequence of random row
vectors $v_{\tau_{i-1}}\dvtx \Omega\to\mathbb{R}^d$
measurable w.r.t. $\mathcal{F}_{\tau_{i-1}} $.
\end{definition}


Let us briefly describe the contents of this paper, which continues and
extends results
from the preprint \cite{Toivola1}.
\begin{longlist}[(2)]
\item[(1)] In Theorem~\ref{thmBesovSpacesNormEquiv}, we provide a
characterization of functions
$f\dvtx \mathbb{R}^d\to\mathbb{R}$
belonging to the Besov space $\mathbb{B}_{p,q}^\theta(\mathbb
{R}^d,\gamma_d)$
by $F\dvtx [0,1]\times\mathbb{R}^d\to\mathbb{R}$ with $F(t,x):= \mathbb
{E}(f(W_1)|W_t=x)$.
Roughly speaking, considering $F$ as the heat
extension of $f\in L_2(\mathbb{R}^d,\gamma_d)$, the regularity of
this extension precisely describes
the Besov regularity of $f$. Theorem~\ref{thmBesovSpacesNormEquiv}
mainly relies on
Proposition~\ref{propgeneralinterpol}, which might be of independent
interest.

\item[(2)]
Besides the real interpolation spaces $\mathbb{B}_{p,q}^\theta$, the
Riemann--Liouville operator $D^{Y,\theta}$ from Section~\ref{secDtheta-operator}
provides an alternative way to describe the fractional regularity of a
function $g\dvtx E\to\mathbb{R}$. It is defined as a functional
of the Hessian matrices $(D^2 G(t,y))_{t\in[0,1)}$ by
\[
D^{Y,\theta}_t g(Y_1):= \biggl( \int
_0^t (1-u)^{1-\theta} H_G^2(u,Y_u)
\,du \biggr) ^{1/2}
\]
with
\[
H_G^2(u,y):= \sum_{k,l=1}^d
\biggl\llvert\biggl(\sigma_{kk}\sigma_{ll}
\frac
{\partial^2 G}{\partial y_k\,\partial y_l} \biggr) (u,y) \biggr\rrvert^2.
\]
In\vspace*{1.5pt} Proposition~\ref{propDthetaAndBesov}, we relate $D^{Y,\theta}$ to
the spaces
$\mathbb{B}_{p,q}^\theta(\mathbb{R}^d,\gamma_d)$, which continues
the analysis in \cite{Gei22}, where
this operator was used in a different form.

\item[(3)]
In the literature
\cite{Gobet-Temam,Gei15,Gei18,MartiniPatry1999,Gei20,Fukasawa2011,Fukasawa2012,Fukasawa2011a}, the question of the behavior of the discretization error
$C_t(g(Y_1),\tau,v)$ has been treated mostly using the \mbox{$L_2$-}norm,
$\| C_t(g(Y_1),\tau,v) \|_{L_2}$, or by weak or stable limits of the
re-scaled error processes
$\lim_n \sqrt{n} C_t(g(Y_1),\tau_n,v_n)$, where $\tau_n$ is of
cardinality $n+1$.
Many of the $L_2$-results are of asymptotic nature as well, and,
concerning random
time-nets, only asymptotic statements were obtained.
There is a general lower bound
$\| C_1(g(Y_1);\tau) \|_{L_p} \ge\delta/\sqrt{n}$ for nets $\tau$
of cardinality $n+1$, see
Remark~\ref{remarklowerbound}. As one of the main results of this
paper, we obtain
in Theorem~\ref{thmLp-adaptedtimenets} a characterization when
this lower bound is actually achieved
by special time-nets. A particular case of this statement is:
\end{longlist}
%

\begin{theorem}\label{prethmLp-adaptedtimenets}
For $2 \le p < \infty$, $0<\theta\le1$, and $g(Y_1)\in L_p$, the
following assertions are equivalent:
\begin{longlist}[(ii)]
\item[(i)] $\| D_1^{Y,\theta}g(Y_1)\|_{L_p} < \infty$,

\item[(ii)] $\sup_{n=1,2,\ldots} \sqrt{n} \| C_1(g(Y_1),\tau_n^\theta) \|
_{L_p} < \infty$,
where the time-nets $\tau_n^\theta=\break (t_{i,n}^\theta)_{i=0}^n$
are given by $t_{i,n}^\theta:= 1 - ( 1-\frac{i}{n} )
^{1/\theta}$.
\end{longlist}
\end{theorem}

Using Proposition~\ref{propDthetaAndBesov}(iii), we can replace, in
the case $p=2$ and $0<\theta<1$, condition (i) in\vspace*{1pt}
the theorem above by $f\in\mathbb{B}_{2,2}^\theta$ [with convention
(\ref{eqnconv2})],
which is in accordance with the known one-dimensional $L_2$-case;
see \cite{Gei19}.

The point of Theorem~\ref{prethmLp-adaptedtimenets} is the usage
of adapted time-nets. In the
literature, equidistant time-nets are often used in discretizations
for simplicity.
Therefore, we provide in Theorem~\ref{thmLp-non-adaptedtimenets} a
description of the random
variables that can be approximated in $L_p$ with equidistant\vspace*{1pt} time-nets
with a rate
$n^{-\theta/2}$ for $0<\theta<1$
in terms of the Besov spaces $\mathbb{B}_{p,\infty}^\theta$. In\vspace*{1pt}
particular, this theorem shows the loss of accuracy
in the approximation
when not using the optimal nets. A special case of this theorem is
the following.

%
\begin{theorem}\label{prethmLp-non-adaptedtimenets}
For $2 \le p < \infty$, $0<\theta<1$ and $g(Y_1)\in L_p$ the
following assertions are equivalent:
\begin{longlist}[(ii)]
\item[(i)] $f \in\mathbb{B}^{\theta}_{p,\infty}$ with $f$ given by (\ref
{eqnconv2}),
\item[(ii)] $\sup_{n=1,2,\ldots} n^{\theta/2} \| C_1(g(Y_1);\tau
_n) \|_{L_p}<\infty$, where
$\tau_n=(i/n)_{i=0}^n$ are the equidistant time-nets.
\end{longlist}
\end{theorem}

(4)~Theorems~\ref{prethmLp-adaptedtimenets} and~\ref{prethmLp-non-adaptedtimenets}
(Theorems~\ref{thmLp-adaptedtimenets} and~\ref{thmLp-non-adaptedtimenets}) are based
on Theorem~\ref{thmLp-error-equivalence} which extends the curvature
type description of
the $L_2$-approximation error from~\cite{Gei15} for deterministic nets
to the $L_p$-error
$\| C_1(g(Y_1),\tau) \|_{L_p}$ with $2\le p< \infty$ and to random time-nets
$\tau=(\tau_i)_{i=0}^n \in{\mathcal T}^{\mathrm{rand}}$. To
illustrate Theorem~\ref{thmLp-error-equivalence}, let us
formulate a corollary that follows from
Remark~\ref{remarkthmLp-error-equivalence}.

%
\begin{theorem}\label{prethmLp-error-equivalence}
For $2 \le p < \infty$ there\vspace*{2pt} is a constant
$c_{\fontsize{8.36pt}{10pt}\selectfont{(\ref{prethmLp-error-equivalence} )}} \ge1$
depending at most on $p$
such that for all $g(Y_1)\in L_p$ and
$\tau=(\tau_i)_{i=0}^n\in{\mathcal T}^{\mathrm{rand}}$ we have that
\[
\bigl\| C_1\bigl(g(Y_1),\tau\bigr) \bigr\|_{L_p}
\sim_{c_{\fontsize{6.6pt}{10pt}\selectfont{(\ref{prethmLp-error-equivalence} )}}} \Biggl\llVert\Biggl(
\sum_{i=1}^n
\int_{\tau_{i-1}}^{\tau_i} (\tau_i-t)
H_G^2(t,Y_t) \,dt \Biggr) ^{1/2} \Biggr\rrVert_{L_p}.
\]
\end{theorem}

For example, to connect Theorem~\ref{prethmLp-error-equivalence} to
Theorem~\ref{prethmLp-adaptedtimenets}, we measure the size of
a sequence $0=t_0 \le\cdots\le t_{n-1}<t_n=1$ by
\[
\bigl|(t_i)_{i=0}^n\bigr|_\theta:=
\frac{|t_{i}-t_{i-1}|}{(1-t_{i-1})^{1-\theta}} \qquad\mbox{with }
0<\theta\le1.
\]
We\vspace*{1.5pt} get $|\tau_n^\theta|_\theta\le1/(\theta n)$, where the nets
$\tau_n^\theta$ are taken from
Theorem~\ref{prethmLp-adaptedtimenets}, in contrast to
$|(i/n)_{i=0}^n|_\theta= n^{-\theta}$
for the equidistant nets, and
Theorem~\ref{prethmLp-error-equivalence} yields
\[
\bigl\| C_1\bigl(g(Y_1),\tau\bigr) \bigr\|_{L_p} \le
c_{\fontsize{8.36pt}{10pt}\selectfont{(\ref{prethmLp-error-equivalence} )}} \bigl\| \sqrt{|\tau|_\theta}
 D_1^{Y,\theta}g(Y_1)
\bigr\|_{L_p} \qquad\mbox{for }\tau\in{\mathcal T}^{\mathrm{rand}}
\]
so that the implication (i)${}\Rightarrow{}$(ii) of Theorem~\ref{prethmLp-adaptedtimenets} follows.

The novelty of Theorem~\ref{prethmLp-error-equivalence}
(Theorem~\ref{thmLp-error-equivalence})
concerns the range $2<p<\infty$ and the fact
that certain fixed random nets (including all deterministic time-nets)
are allowed, which distinguishes
the result from previous asymptotic ones.
As already pointed out, the techniques for the
\mbox{$L_p$-}estimates differ significantly from the \mbox{$L_2$-}estimates because
the problem cannot be
translated into a one-step approximation problem nor can we use orthogonality.
Moreover, the extension from deterministic nets
to random nets does not seem to be straightforward as we still have to
use the sub-class
${\mathcal T}^{\mathrm{rand}}$ of random nets $(\tau_i)_{i=1}^n$
where $\tau_i$ is
$\mathcal{F}_{\tau_{i-1}}$-measurable
(see Remark~\ref{remarkrestrictiontimenets}).
Our \mbox{$L_p$-}estimates can be seen as an interpolation between the
$L_2$-estimates mentioned above
and the weighted BMO-estimates from \cite{Gei17}.
However, pure interpolation techniques do not seem to be sufficient yet
to fully treat our problem.

\section{Preliminaries}\label{secprelimnaries}
\subsection*{Notation}
We use $A \sim_c B$ for $A/c \le B \le cA$ whenever $A,B\ge0$ and
\mbox{$c\ge1$},
$a\vee b= \max\{a,b\}$ and $a\wedge b=\min\{a,b\}$, and
let $|\cdot|$ be the Euclidean norm for a vector or the
Hilbert--Schmidt norm for a
matrix. Given a random vector or a random matrix $A$, we write
$\|A\|_{L_p}:= \||A|\|_{L_p}$ and denote the transpose of $A$ by~$A^\top$.

\subsection*{Real interpolation}
Let us recall the real interpolation method
that we use to generate the (Gaussian) Besov spaces.

%
\begin{definition}[(\cite{B-S,B-L})]
Let $(X_0,X_1)$ be a compatible couple of Banach spaces,
that is, there exists a Hausdorff topological vector space in which
both $X_0$ and $X_1$ are continuously embedded.
Given $x\in X_0+X_1$ and $\lambda> 0$, the \mbox{$K$-}functional is defined by
\[
K(x,\lambda;X_0,X_1):= \inf\bigl\{ \| x_0
\|_{X_0} + \lambda\| x_1\|_{X_1}\dvtx
x=x_0+x_1, x_i\in X_i \bigr\}.
\]
Given $0<\theta<1$ and $1\le q \le\infty$, we let
$(X_0,X_1)_{\theta,q}$ be the space of all $x\in X_0+X_1$ such that
\[
\| x \|_{(X_0,X_1)_{\theta,q}}:= \bigl\llVert\lambda^{-\theta} K(x,
\lambda;X_0,X_1) \bigr\rrVert_{L_q((0,\infty),(d \lambda)/\lambda
)} < \infty.
\]
\end{definition}

The $K$-functional yields to one of the basic approaches to define
intermediate spaces $Y$ of a compatible couple of Banach spaces
$(X_0,X_1)$, that is, Banach spaces $Y$ such that one has continues
embeddings $X_0\cap X_1 \hookrightarrow Y \hookrightarrow X_0+X_1$.
Assuming $X_1 \hookrightarrow X_0$ with
norm one, this reduces to the embedding $X_1 \hookrightarrow Y
\hookrightarrow X_0$. In this case,
$K(x,\lambda;X_0,X_1)= \|x\|_{X_0}$ for $\lambda\in[1,\infty)$
which does not give any information. However, for
$\lambda\in(0,1)$ we have that
\[
\lambda\|x\|_{X_0} \le K(x,\lambda;X_0,X_1)
\le\|x\|_{X_0}.
\]
The behavior of the function $\lambda\to K(x,\lambda;X_0,X_1)$ close
to zero describes the distance of
$x$ to $X_1$: intuitively we can say that the closer the function is to
a linear function in $\lambda$,
the closer $x$ is to $X_1$. In general, without the restriction $X_1
\hookrightarrow X_0$,
the functionals
\[
\bigl\llVert\lambda^{-\theta} K(x,\lambda;X_0,X_1)
\bigr\rrVert_{L_q((0,\infty),(d \lambda)/\lambda)}
\]
examine the behavior of the $K$-functional (in particular at zero and
at infinity)
and lead to the spaces
$(X_0,X_1)_{\theta,q}$. For $X_1 \hookrightarrow X_0$, we obtain the
lexicographical ordering
\[
(X_0,X_1)_{\theta_0,q_0} \subseteq(X_0,X_1)_{\theta_1,q_1}
\quad\mbox{and}\quad(X_0,X_1)_{\eta,r_0}
\subseteq(X_0,X_1)_{\eta,r_1}
\]
if $0<\theta_1<\theta_0<1$, $1\le q_0,q_1 \le\infty$,
$0<\eta<1$, and $1\le r_0\le r_1 \le\infty$.
The choice of the
measure $d\lambda/ \lambda$ ensures (also in the general case) the symmetry
$(X_0,X_1)_{\theta,q} = (X_1,X_0)_{1-\theta,q}$.

\subsection*{Gaussian Sobolev and Besov spaces}
We let $d\ge1$ and
$\gamma_d$ be the standard Gaussian measure on $\mathbb{R}^d$. The
space $L_2(\mathbb{R}^d,\gamma_d)$ is equipped
with the orthonormal basis of generalized Hermite polynomials
$(h_{k_1,\ldots,k_d})_{k_1,\ldots,k_d=0}^\infty$
given by
\[
h_{k_1,\ldots,k_d}(x_1,\ldots,x_d):=
h_{k_1}(x_1)\cdots h_{k_d}(x_d),
\]
where $(h_k)_{k=0}^\infty\subset L_2(\mathbb{R},\gamma_1)$ is the
standard orthonormal basis of Hermite polynomials.
The Sobolev space $\mathbb{D}_{1,2} = \mathbb{D}_{1,2}(\mathbb
{R}^d,\gamma_d)$ consists of all $f\in L_2 (\mathbb{R}^d,\break \gamma_d)$
such that
\[
\sum_{k_1,\ldots,k_d=0}^\infty\langle f,
h_{k_1,\ldots,k_d} \rangle_{L_2(\mathbb{R}^d,\gamma_d)}^2 \llVert
\nabla
h_{k_1,\ldots,k_d} \rrVert_{L_2(\mathbb{R}^d,\gamma_d)}^2 < \infty.
\]
The space $\mathbb{D}_{1,2}$ is a Banach space under the norm
\[
\| f \|_{\mathbb{D}_{1,2}}:= \sqrt{\| f\|_{L_2(\mathbb{R}^d,\gamma
_d)}^2 +
\| Df \|_{L_2(\mathbb{R}^d,\gamma_d)}^2},
\]
where, for $f \in\mathbb{D}_{1,2}$, the gradient $Df$ is given by
\[
D f:= \sum_{k_1,\ldots,k_d=0}^\infty\langle f,
h_{k_1,\ldots,k_d} \rangle_{L_2(\mathbb{R}^d,\gamma_d)} \nabla
h_{k_1,\ldots,k_d}.
\]
Given $2\le p < \infty$, the Banach space $\mathbb{D}_{1,p} \subseteq
L_p$ is given by
\[
\mathbb{D}_{1,p}:= \bigl\{ f \in\mathbb{D}_{1,2}\dvtx  \| f
\|_{\mathbb{D}_{1,p}}:= \bigl( \| f \|_{L_p}^p + \| Df \|
_{L_p}^p \bigr) ^{1/p} < \infty\bigr\}.
\]
Here and later, we use $\|f\|_{L_p} = \| f\|_{L_p(\mathbb{R}^d,\gamma
_d)}$ and
$\|Df\|_{L_p} = \| Df\|_{L_p(\mathbb{R}^d,\gamma_d)}$.

%
\begin{definition}
For $0<\theta< 1$ and $1\le q \le\infty$, we let
\[
\mathbb{B}_{p,q}^\theta:= (L_p,
\mathbb{D}_{1,p})_{\theta,q}
\]
be the Gaussian Besov space on $\mathbb{R}^d$ of fractional smoothness
$\theta$ and
fine-index $q$.
\end{definition}

Because $\mathbb{D}_{1,p}$ is not closed in $L_p$, we get a scale of
spaces indexed by $(\theta,q)$,
where the spaces are identical if and only if both indices coincide
(see \cite{Jansonetal1984}, Theorem 3.1). A typical function
which has fractional smoothness is given by the following.

%
\begin{example}\label{exampleBesovspaces}
Let $d=1$, $K\in\mathbb{R}$, $2\le p < \infty$ and $0 \le\alpha<
1- \frac{1}{p}$. Then one has that
\[
f(x):= \cases{ \bigl((x-K)^+\bigr)^\alpha, &\quad $\alpha> 0$
\vspace*{2pt}\cr
\chi_{[K,\infty)}(x), &\quad  $\alpha= 0$}\in\mathbb
{B}_{p,\infty}^{(1/p)+\alpha},
\]
which shows the trade-off between integrability
and smoothness. This can be proved by verifying
\[
K(f,\lambda;L_p,\mathbb{D}_{1,p}) \le c
\lambda^{(1/p)+\alpha} \qquad\mbox{for }0<\lambda<1.
\]
\end{example}

Using canonical representations of functions of bounded variation,
one can extend the case $\alpha=0$ in Example \ref
{exampleBesovspaces} to certain
functions of bounded variation by considering convex combinations
$f(x)=\sum_{l=1}^L \beta_l \chi_{[K_l,\infty)}(x)$.

\subsection*{Burkholder--Davis--Gundy inequality}
We use the Burkholder--Davis--Gun\-dy inequality for Brownian
martingales with
values in a separable Hilbert space. An explicit formulation is as follows:
assume for $i=1,2,\ldots$ progressively measurable processes
$(L_t^i)_{t\in[0,1]}$ with
$L_t^i\dvtx \Omega\to\mathbb{R}^d$ considered as row vectors and such that
\[
\sum_{i=1}^\infty\mathbb{E}\int
_0^1 \bigl|L_t^i\bigr|^2
\,dt < \infty,
\]
then, for all $1<p<\infty$, there is a constant $c_{\fontsize{8.36pt}{10pt}\selectfont{(\ref{eqnBDG})}}=
c_{\fontsize{8.36pt}{10pt}\selectfont{(\ref{eqnBDG})}}(p)\ge1$ such that
%
\begin{equation}
\label{eqnBDG} \Biggl\|\Biggl( \sum
_{i=1}^\infty\biggl\llvert\int_0^1
L_u^i \,dW_u \biggr\rrvert^2
\Biggr) ^{1/2} \Biggr\|_{L_{p}}
\sim_{c_{\fontsize{6.6pt}{10pt}\selectfont{(\ref{eqnBDG})}}} \Biggl\|\Biggl( \sum
_{i=1}^\infty\int_0^1
\bigl|L_u^i\bigr|^2 \,du \Biggr)^{1/2} \Biggr\Vert_{L_{p}}.
\end{equation}


\section{Fractional smoothness on the Gaussian space}\label{secfractionalsmoothnesswienerspace}

In this section, we characterize the Gaussian Besov spaces $\mathbb
{B}^{\theta}_{p,q}$
by the behavior of $G$ from (\ref{eqndefinitionG}) in the case
$Y=W$. To make
this more clear, we do a change of notation and replace $g$ by $f$ and
$G$ by $F$.
This means that $f\in L_2(\mathbb{R}^d,\gamma_d)$ and
\[
F(t,x):= \mathbb{E}f(x+W_{1-t}) \qquad\mbox{for }(t,x)\in[0,1]
\times\mathbb{R}^d.
\]
We also use the Hessian $d\times d$ matrix
\[
D^2F:= \biggl( \frac{\partial^2 F}{\partial x_i\,\partial x_j} \biggr
)_{i,j=1}^d.
\]
One can check that
%
\begin{equation}
\label{eqnd12hessian} f\in\mathbb{D}_{1,2} \qquad\mbox{if and only
if }\int_0^1 \bigl\| D^2
F(t,W_t)\bigr\|_{L_2}^2 \,dt < \infty.
\end{equation}
Moreover, for all $f\in L_2(\mathbb{R}^d,\gamma_d)$ we have
%
\begin{equation}
\nabla F(t,W_t) = \nabla F(0,0) + \biggl( \int_0^t
D^2 F(u,W_u) \,dW_u \biggr) ^\top
\label{eqnrepresentationnablaF} \qquad\mbox{a.s.}
\end{equation}
for $0\le t < 1$, where $\nabla F(t,x)$ is considered as a row vector.
If $f\in\mathbb{D}_{1,2}$, then (\ref{eqnrepresentationnablaF})
can be extended to $t=1$ with the convention
$\nabla F(1,\cdot):= Df$.
Now we generalize~(\ref{eqnd12hessian}) to the scale of Besov spaces.

%
\begin{theorem} \label{thmBesovSpacesNormEquiv}
Let $2 \leq p < \infty$, $0 < \theta< 1$, $1 \leq q \leq\infty
$
and $f\in L_p(\mathbb{R}^d,\gamma_d)$. Then
\begin{eqnarray*}
\|f \|_{\mathbb{B}^{\theta
}_{p,q}} &\sim_{c_{\fontsize{6.6pt}{10pt}\selectfont{(\ref{thmBesovSpacesNormEquiv})}}} & \|f \|
_{L_{p}} + \bigl\|(1-t)^{-\theta/2} \bigl\Vert F(1, W_1) - F(t, W_t) \bigr\Vert_{L_{p}} \bigr\|_{L_{q} ([0,1), (dt)/(1-t) )}
\\
& \sim_{c_{\fontsize{6.6pt}{10pt}\selectfont{(\ref{thmBesovSpacesNormEquiv})}}}& \|f \|
_{L_{p}} + \bigl\|(1-t)^{(1-\theta)/2} \bigl\Vert\nabla F(t, W_t) \bigr\|
_{L_{p}} \bigr\|_{L_{q} ([0,1), (dt)/(1-t))}
\\
& \sim_{c_{\fontsize{6.6pt}{10pt}\selectfont{(\ref{thmBesovSpacesNormEquiv})}}} & \|f \|
_{L_{p}} + \bigl\|(1-t)^{(2-\theta)/2} \bigl\Vert D^2 F(t, W_t) \bigr\|
_{L_{p}} \bigr\|_{L_{q} ([0,1), (dt)/(1-t))},
\end{eqnarray*}
where $c_{\fontsize{8.36pt}{10pt}\selectfont{(\ref{thmBesovSpacesNormEquiv})}} \geq1$ depends uniquely
on $(p,\theta,q)$.
\end{theorem}

%
\begin{remark}
Theorem~\ref{thmBesovSpacesNormEquiv} generalizes \cite{Gei19}, Theorem
2.2, where $p=2$ was considered,
and \cite{Toivola1}, Lemma 4.7, which was proved for $2<p<\infty$ and
$q=\infty$.
\end{remark}

Before we prove Theorem~\ref{thmBesovSpacesNormEquiv}, we derive
a corollary in the case $d=1$ concerning the oscillation of a Borel function
$f\dvtx \mathbb{R}\to\mathbb{R}$ given by
\[
\mathrm{OSC}_p(f,x_0,s):= \biggl( \frac{1}{4s^2} \int
_{Q(x_0,s)} \bigl|f(y)-f(z)\bigr|^p \,dy \,dz \biggr)^{1/p},
\]
where $2\le p< \infty$, $s>0$, $x_0\in\mathbb{R}$ and
$Q(x_0,s):= \{ (y,z)\dvtx  |y-x_0| \le s, |z-x_0| \le s \}$.

\begin{corollary}\label{corosc}
For $2\le p < \infty$, $0<\theta<1$, $1\le q \le\infty$ and $f\in
\mathbb{B}^{\theta}_{p,q}$,
we have that
\[
\bigl\llVert s^{(1/p)-\theta} \mathrm{OSC}_p(f,x_0,s)
\bigr\rrVert_{L_q((0,1],(ds)/s)} \le c_{\fontsize{8.36pt}{10pt}\selectfont{(\ref{corosc})}}
\| f \|_{\mathbb
{B}^{\theta}_{p,q}},
\]
where the constant $c_{\fontsize{8.36pt}{10pt}\selectfont{(\ref{corosc})}}>0$ depends at most on
$(p,\theta,q,x_0)$.
\end{corollary}

\begin{pf}
From \cite{Gei22}, Lemma 4.9, we know that
\[
\mathrm{OSC}_p(f,x_0,\sqrt{1-t}) \le c
(1-t)^{-1/(2p)} \bigl\| f(Y)-f(Z) \bigr\|_{L_p}
\]
for $f\in L_p(\mathbb{R},\gamma_1)$, $0\le t < 1$ and
a two-dimensional Gaussian vector $(Y,Z)$ with $Y,Z\sim N(0,1)$ and
$\operatorname{cov}(Y,Z) = t$, where $c>0$ depends at most on $(x_0,p)$. Looking at
\cite{Gei22}, Proof of Proposition 4.5(iii), we see that
\[
\bigl\| f(Y)-f(Z) \bigr\|_{L_p} \le2 \bigl\| f(W_1) - \mathbb{E}
\bigl(f(W_1) | \mathcal{F}_t \bigr) \bigr\|_{L_p},
\]
so that we can conclude by Theorem~\ref{thmBesovSpacesNormEquiv} of
this paper.
\end{pf}
Now we turn to the proof of Theorem~\ref{thmBesovSpacesNormEquiv}. We
start with the following proposition.

%
\begin{proposition} \label{propestimateKfunctionalNoDerivative}
Let $2\le p < \infty$. There exists a constant
$ c_{\fontsize{8.36pt}{10pt}\selectfont{(\ref{propestimateKfunctionalNoDerivative} )}} \geq
1 $ depending at most on $p$ such that
for any $0<t<1$,
\[
K(f, \sqrt{1-t}; L_p, \mathbb{D}_{1,p})
\sim_{c_{\fontsize{6.6pt}{10pt}\selectfont{(\ref{propestimateKfunctionalNoDerivative} )}}} \bigl( \bigl\| f(W_1) - F(t,
W_t) \bigr\|_{L_{p}} + \sqrt{1-t} \| f \|_{L_{p}} \bigr).
\]
\end{proposition}

\begin{pf}
(a) Fix $0<t<1$ and $\varepsilon>0$. We find
$f_0 \in L_p$ and $f_1 \in\mathbb{D}_{1,p} $
such that $f = f_0 + f_1$ and
\[
\|f_0 \|_{L_p} +
\sqrt{1-t} \|f_1 \|
_{\mathbb{D}_{1,p}} \leq K(f,\sqrt{1-t};L_p,\mathbb{D}_{1,p})
+ \epsilon.
\]
For $F_i(t,x):= \mathbb{E}(f_i(W_1)|W_t=x)$ we obtain from (\ref
{eqnItoG}) and (\ref{eqnBDG}) that
\begin{eqnarray*}
&& \bigl\| f(W_1) - F(t, W_t) \bigr\|_{L_{p}}
\\
&&\qquad \leq \bigl\| f_0(W_1) -
F_0(t, W_t) \bigr\|_{L_{p}} +
\biggl\|\int_t^1
\nabla{F_1}(u, W_u) \,dW_u \biggr\|_{L_{p}}
\\
&&\qquad \leq \bigl\| f_0(W_1) -
F_0(t, W_t) \bigr\|_{L_{p}} +
c_{\fontsize{8.36pt}{10pt}\selectfont{(\ref{eqnBDG})}} \biggl( \int_t^1 \bigl\|
\nabla{F_1}(u, W_u) \bigr\Vert_{L_{p}} ^2 \,du \biggr)^{1/2}
\\
&&\qquad \leq 2 \|f_0\|_{L_{p}} + c_{\fontsize{8.36pt}{10pt}\selectfont{(\ref{eqnBDG})}} \sqrt{1-t} \|f_1\|_{\mathbb{D}_{1,p} }
\\
&&\qquad \leq c \bigl[ K(f, \sqrt{1-t};L_p,\mathbb{D}_{1,p}) +
\epsilon\bigr],
\end{eqnarray*}
where $c:= \max\{c_{\fontsize{8.36pt}{10pt}\selectfont{(\ref{eqnBDG})}}, 2\}$ and we employed the facts
that $2 \leq p < \infty$ and that
(\ref{eqnrepresentationnablaF}) yields
\[
\bigl\|\nabla{F_1}(u, W_u)\bigr\|_{L_{p}} \leq\|f_1\|_{\mathbb{D}_{1,p} }\qquad\mbox{for all } 0 \leq
u \leq1.
\]
Letting $\epsilon\to0$
and observing that $\sqrt{1-t} \| f \|
_{L_{p}} \leq K(f, \sqrt{1-t}; L_p, \mathbb{D}_{1,p})$
we achieve the first part of the desired inequality.

(b) For $0<t<1$, we set
\[
g_t (x):= F(t, \sqrt{t} x) \quad\mbox{and}\quad h_t (x):= f(x) - F(t,
\sqrt{t} x)
\]
so that
\begin{eqnarray*}
\| g_t \|_{\mathbb{D}_{1,p} }^p
& = & \bigl\| F(t, \sqrt{t} x) \bigr\Vert
_{L_p}^p + \bigl\|\nabla F(t, \sqrt{t}
x) \sqrt{t} \bigr\|_{L_p}^p
\\
&\le& \|f \|_{L_p}^p
+ \bigl\|\nabla F(t, W_t) \bigr\Vert_{L_p}^p.
\end{eqnarray*}
Applying (\ref{eqnItoG}) for $Y=W$, (\ref{eqnBDG}), the fact that
$\|\nabla F(t, W_t) \|_{L_{p}}$ is
nondecreasing in $t$ and that $2\le p <\infty$,
we estimate
\begin{eqnarray*}
\bigl\| F(t, W_t) \bigr\|
_{L_{p}} & \leq& \biggl\|\int_0^t
\nabla F(u, W_u) \,dW_u \biggr\|
_{L_{p}} + \bigl\| F(0, W_0) \bigr\|_{L_{p}}
\\
& \leq& c_{\fontsize{8.36pt}{10pt}\selectfont{(\ref{eqnBDG})}} \bigl\|\nabla F(t, W_t)
\bigr\|_{L_{p}} + \bigl| \mathbb{E}f(W_1) \bigr|.
\end{eqnarray*}
Thus,
\begin{eqnarray*}
\| g_t \|_{\mathbb{D}_{1,p} }
&\le& \bigl\| f(W_1) - F(t, W_t) \bigr\|_{L_{p}} + (1 + c_{\fontsize{8.36pt}{10pt}\selectfont{(\ref{eqnBDG})}}) \bigl\|\nabla F(t, W_t) \bigr\|_{L_{p}} +
\bigl| \mathbb{E}f(W_1) \bigr|
\\
&\le& \bigl[ 1 + (1+c_{\fontsize{8.36pt}{10pt}\selectfont{(\ref{eqnBDG})}}) c_{\fontsize{8.36pt}{10pt}\selectfont{(\ref{lemmaestimateDerFbyf} )}} (1-t)^{-1/2}
\bigr] \bigl\| f(W_1) - F(t, W_t) \bigr\|_{L_{p}} + \bigl| \mathbb{E}f(W_1) \bigr|,
\end{eqnarray*}
where we used Lemma~\ref{lemmaestimateDerFbyf}.
Exploiting an independent Brownian motion $\widetilde{W}$ and the fact
that the covariance structures of $(W_1,\sqrt{t} W_1 + \sqrt{1-t}\widetilde{W}_1)$ and
$(W_1, W_{\sqrt{t}} + \widetilde{W}_{1-\sqrt{t}})$ are the same, we
obtain for $h_t$
that
\begin{eqnarray*}
\| h_t \|_{L_p} & = &
\bigl[ \mathbb{E}\bigl\llvert f(W_1) - \widetilde{\mathbb{E}} f(\sqrt{t}
W_1 + \sqrt{1-t}\widetilde{W}_1) \bigr\rrvert^p
\bigr]^{1/p}
\\
&\leq& \bigl[ \mathbb{E}\widetilde{\mathbb{E}} \bigl\llvert f(W_1) - f(
W_{\sqrt{t}} + \widetilde{W}_{1-\sqrt{t}}) \bigr\rrvert^p
\bigr]^{1/p}
\\
&\leq& \bigl\| f(W_1) - F(\sqrt{t},
W_{\sqrt{t}}) \bigr\|_{L_{p}} + \bigl\| F(\sqrt{t}, W_{\sqrt{t}}) - f( W_{\sqrt{t}} +
\widetilde{W}_{1-\sqrt{t}}) \bigr\|_{L_{p}}
\\
& = & 2 \bigl\| f(W_1) - F(\sqrt{t},
W_{\sqrt{t}}) \bigr\|_{L_{p}}
\\
&\leq& 4 \bigl\| f(W_1) - F(t, W_t)
\bigr\|_{L_{p}},
\end{eqnarray*}
where in the last step $F(t, W_t)$ was inserted. Hence,
\begin{eqnarray*}
&& K(f, \sqrt{1-t}; L_p, \mathbb{D}_{1,p})
\\
&&\qquad \le \| h_t \|_{L_p} +
(1-t)^{1/2} \| g_t\|_{\mathbb
{D}_{1,p}}
\\
&&\qquad \le (1-t)^{1/2} \bigl|\mathbb{E}f(W_1)\bigr|
 + \bigl[5+(1+c_{\fontsize{8.36pt}{10pt}\selectfont{(\ref{eqnBDG})}})c_{\fontsize{8.36pt}{10pt}\selectfont{(\ref{lemmaestimateDerFbyf} )}}\bigr] \bigl\| f(W_1) - F(t, W_t) \bigr\|
_{L_{p}}
\end{eqnarray*}
and the proof is complete.
\end{pf}

We are now ready to prove the main result of this section.

\begin{pf*}{Proof of Theorem~\ref{thmBesovSpacesNormEquiv}}
To verify the assumptions of Proposition~\ref{propgeneralinterpol},
we set
\begin{eqnarray*}
d^0(t) &:= & \bigl\| f(W_1) - F(t,W_t)
\bigr\|_{L_p},
\\
d^1(t) &:= & \bigl\| \nabla F(t,W_t) \bigr\|_{L_p},
\\
d^2(t) &:= & \bigl\| D^2 F(t,W_t)
\bigr\|_{L_p},
\end{eqnarray*}
$A:= 2 c_{\fontsize{8.36pt}{10pt}\selectfont{(\ref{lemmaestimateDerFbyf} )}}
\| f \|_{L_{p}}$ and
$\alpha:= c_{\fontsize{8.36pt}{10pt}\selectfont{(\ref{eqnBDG})}} \vee c_{\fontsize{8.36pt}{10pt}\selectfont{(\ref{lemmaestimateDerFbyf})}}$.
Then Lemma~\ref{lemmaestimateDerFbyf} implies that
\[
d^k(t) \le c_{\fontsize{8.36pt}{10pt}\selectfont{(\ref{lemmaestimateDerFbyf})}} (1-t)^{-k/2}
\,d^0(t) \qquad\mbox{for } k=1,2.
\]
By (\ref{eqnItoG}), (\ref{eqnrepresentationnablaF}),
the Burkholder--Davis--Gundy inequalities (\ref{eqnBDG})
and $2\le p < \infty$, we also see that
\[
d^0(t) \le c_{\fontsize{8.36pt}{10pt}\selectfont{(\ref{eqnBDG})}} \biggl( \int_0^t
\bigl\llVert\nabla F(s,W_s) \bigr\rrVert_{L_p}^2
\,ds \biggr)^{1/2} = c_{\fontsize{8.36pt}{10pt}\selectfont{(\ref{eqnBDG})}} \biggl( \int_0^t
\bigl[d^1(s)\bigr]^2 \,ds \biggr)^{1/2}
\]
and
\begin{eqnarray*}
d^1(t) &\le& \bigl\llVert\nabla F (0, W_0) \bigr
\rrVert_{L_p} + \biggl\llVert\int_0^t
D^2 F(s,W_s) \,dW_s \biggr\rrVert
_{L_p}
\\
&\le& 2 c_{\fontsize{8.36pt}{10pt}\selectfont{(\ref{lemmaestimateDerFbyf} )}} \|f \|
_{L_{p}} + c_{\fontsize{8.36pt}{10pt}\selectfont{(\ref{eqnBDG})}} \biggl( \int_0^t
\bigl\llVert D^2 F(s,W_s) \bigr\rrVert
_{L_p}^2 \,ds \biggr)^{1/2}
\\
& = & 2 c_{\fontsize{8.36pt}{10pt}\selectfont{(\ref{lemmaestimateDerFbyf} )}} \|f \|
_{L_{p}} + c_{\fontsize{8.36pt}{10pt}\selectfont{(\ref{eqnBDG})}} \biggl( \int_0^t
\bigl[d^2(s)\bigr]^2 \,ds \biggr)^{1/2},
\end{eqnarray*}
where we used Lemma~\ref{lemmaestimateDerFbyf}.
Now, applying (\ref{eqngeneralinterpolcor}) on page \pageref{eqngeneralinterpolcor}
gives the equivalence
between the last three expressions in Theorem \ref
{thmBesovSpacesNormEquiv}. It remains to check that
\[
\|f \|_{\mathbb{B}^{\theta
}_{p,q}}
\sim_c \|f \|_{L_{p}}
+ \bigl\|(1-t)^{-\theta/2} \bigl\| F(1,
W_1) - F(t, W_t) \bigr\|
_{L_{p}} \bigr\|_{L_{q} ([0,1), (dt)/(1-t))}
\]
for some $c=c(p,q,\theta)\ge1$,
which follows from Proposition~\ref{propestimateKfunctionalNoDerivative}.
\end{pf*}

%
\begin{remark}
\label{remarkOU-semigroup}
In the literature, interpolation spaces on the Wiener (or Gaussian)
space are considered,
for example, in \cite{Watanabe,Hirsch,Gei19,Gei22,Toivola1}.
A classical approach is based on semi-groups. Instead of that, our
approach uses the elementary
Proposition~\ref{propgeneralinterpol}, which is not related to
semi-groups and makes it
therefore possible to apply Proposition~\ref{propgeneralinterpol} in
more general situations
(see \cite{geigob2012}).
Regarding the present paper, Proposition~\ref{propgeneralinterpol}
opens the way to extend results from Sections~\ref{secDtheta-operator}
and~\ref{secappoximation-p} below to processes different from the
(geometric) Brownian motion.
Below we want to indicate a possible semi-group approach to Theorem
\ref{thmBesovSpacesNormEquiv}:
\begin{longlist}[(1)]
\item[(1)]
The first equivalence of Theorem~\ref{thmBesovSpacesNormEquiv} can be
deduced in the
case $q=p$ ($q$~is the fine-tuning index in the interpolation, $L_p$
the integrability
of the underlying spaces)
from \cite{Hirsch}, Remark on page~428. Using the simple observation
\cite{Geiss-Geiss-Gobet2012}, equation~(6),
one can transform Hirsch's condition into
\[
\int_0^\infty s^{-(\theta p)/2} \bigl\|
f(W_1) - F\bigl(e^{-s}, W_{e^{-s}}\bigr)
\bigr\|_p^p \frac{ds}{s},
\]
which is our condition, up to a different scaling.

\item[(2)]
To consider the general case, that is, $q\neq p$, and also the other
equivalences in
Theorem~\ref{thmBesovSpacesNormEquiv}, one can check general results
about interpolation
and semi-groups. There are two natural semi-groups one might use,
the Ornstein--Uhlenbeck semi-group and the Poisson semi-group (see
\cite{ste70}).
Roughly speaking, switching from the Ornstein--Uhlenbeck semi-group to the
Poisson semi-group should result in a change of the main interpolation
parameter $\eta$
to our parameter $\theta=\eta/2$ in the corresponding formulas,
cf. \cite{Triebel-interpolation}, Section~1.15.2.
Now assume $2\le p < \infty$ and $\xi=f(W_1) \in L_p$ and let
$(T_s)_{s\ge0}$ be the
Ornstein--Uhlenbeck semi-group on $L_p$ with generator $\Lambda$. Then
\cite{Triebel-interpolation}, Section~1.13.2, gives
%
\begin{equation}
\label{eqnsemigroup-interpolation} \| f\|_p + \bigl\llVert s^{-\eta} \|
T_s \xi- \xi\|_p \bigr\rrVert_{L_q((0,\infty),(ds)/s)} < \infty
\end{equation}
for the interpolation space $(L_p,D(\Lambda))_{\eta,q}$ with $0<\eta
<1$ and $1\le q \le\infty$.
By Mehler's formula, we have
\[
T_s \xi= \widetilde{\mathbb{E}} f\bigl(fe^{-s}
W_1 + \sqrt{1-e^{-2s}} \widetilde W_1\bigr) = F
\bigl(e^{-2s}, e^{-s} W_1\bigr)
\]
for an independent Brownian motion $\widetilde W$. This
would give a comparable statement to the first equivalence
of Theorem~\ref{thmBesovSpacesNormEquiv} for the Ornstein--Uhlenbeck
semi-group.
To come closer to our statement, one can inspect the proof of
Proposition~\ref{propestimateKfunctionalNoDerivative}, which gives
$\| T_s \xi- \xi\|_p \le4 \| f(W_1) - F(e^{-2s},W_{e^{-2s}})\|_p$.
Inserting this upper bound into~(\ref{eqnsemigroup-interpolation})
would give an expression like in the first equivalence
of Theorem~\ref{thmBesovSpacesNormEquiv}.
To get the full statement one would still need to try to upper bound
$\| f(W_1) - F(t,W_t)\|_p$ by $\| T_s \xi- \xi\|_p$ in an appropriate way
or to find an alternative way. Concerning the second and third
equivalence of
Theorem~\ref{thmBesovSpacesNormEquiv} one might try to exploit
\cite{Triebel-interpolation}, Section~1.14.5.
\end{longlist}
\end{remark}


\section{The Riemann--Liouville operator \texorpdfstring{$D^{Y,\theta}$}{DY,theta}}\label{secDtheta-operator}

Riemann--Liouville type operators are typically used
to describe fractional re\-gu\-larity. We use these operators to
replace the Besov regularity defined by real interpolation when we consider
the approximation along adapted time-nets in Theorem \ref
{thmLp-adaptedtimenets}
below.
The operator, introduced in the following Definition~\ref{definitonDbeta},
was also used in a slightly modified form in~\cite{Gei22}, where the weak
convergence of the error processes was considered.

%
\begin{definition}\label{definitonDbeta}
For $g(Y_1)\in L_2$, $0<\theta\le1$ and $0\le t\le1$, we let
\[
D^{Y,\theta}_t g(Y_1):= \biggl( \int
_0^t (1-u)^{1-\theta} H_G^2(u,Y_u)
\,du \biggr) ^{1/2},
\]
where, with $G$ given by (\ref{eqndefinitionG}),
\[
H_G^2(u,y):= \sum_{k,l=1}^d
\biggl\llvert\biggl(\sigma_{kk}\sigma_{ll}
\frac
{\partial^2 G}{\partial y_k\,\partial y_l} \biggr) (u,y) \biggr\rrvert^2.
\]
\end{definition}

From now on, we use the following convention:
for $x=(x_1,\ldots,x_d)^\top\in\mathbb{R}^d$ and
$0\le t \le1$ we let
%
\begin{equation}
\label{eqnconv1} y_k(t) = \cases{ \displaystyle x_k, &\quad $Y =W$
\vspace*{3pt}\cr
\displaystyle e^{x_k-(t/2)}, &\quad else} \quad\mbox{and}\quad
y(t):=
\bigl(y_1(t),\ldots,y_d(t)\bigr)^\top
\end{equation}
and define the functions $f\dvtx \mathbb{R}^d\to\mathbb{R}$ and
$F\dvtx [0,1]\times\mathbb{R}^d\to\mathbb{R}$
as
%
\begin{equation}
\label{eqnconv2} f(x):= g\bigl(y(1)\bigr) \quad\mbox{and}\quad
F(t,x):=
\mathbb{E}f(x+W_{1-t})
\end{equation}
so that $ f(W_1) = g(Y_1)$ and $F(t, x)=G(t, y(t))$.
In the case that $Y$ is the coordinate-wise geometric Brownian
motion, this notation implies that
%
\begin{equation}
\label{eqnDerBMgBM} y_k(t) y_l(t) \frac{\partial^2 G}{\partial
y_k\,\partial y_l}
\bigl(t,y(t)\bigr) = \frac{\partial^2 F}{\partial x_k\,\partial
x_l}(t,x) - \delta_{k,l}
\frac{\partial F}{\partial x_k}(t,x)
\end{equation}
for $k,l=1,\ldots,d$.
Let us summarize the connections between the Besov spaces and the
operator $D^{Y,\theta}$ known to us.

%
\begin{proposition} \label{propDthetaAndBesov}
For $g(Y_1)\in L_p$ with $2\le p < \infty$, the following assertions
hold true:
\begin{longlist}[(iii)]
\item[(i)] If $2 < p < \infty$ and $0 < \theta< 1$, then
\begin{enumerate}[(a)]
\item[(a)]$f\in\mathbb{B}_{p,2}^\theta$ implies
$D^{Y,\theta}_1 g(Y_1)\in L_p$,
\item[(b)]$D^{Y,\theta}_1 g(Y_1)\in L_p$ implies $f\in\mathbb
{B}_{p,\infty}^\theta$.
\end{enumerate}
\item[(ii)] If $2\le p < \infty$, then
$D^{Y,1}_1 g(Y_1)\in L_p$ if and only if $f\in\mathbb{D}_{1,p}$.
\item[(iii)] If $0<\theta<1$, then
$D^{Y,\theta}_1 g(Y_1) \in L_2$ if and only if $f\in\mathbb
{B}_{2,2}^\theta$.
\end{longlist}
\end{proposition}

\begin{pf}
(i)
\begin{longlist}[(a)]
\item[(a)] Because $2\le p < \infty$, we see that
\begin{eqnarray*}
\bigl\| D^{Y,\theta}_1 g(Y_1)
\bigr\|_{L_p} &\le& \biggl( \int_0^1
(1-t)^{1-\theta} \bigl\| H_G(t,Y_t)
\bigr\|_{L_p}^2 \,dt \biggr) ^{1/2}
\\
& = & \bigl\|(1-t)^{(2-\theta)/2}
\bigl\|
H_G(t,Y_t) \bigr\|_{L_{p}}
\bigr\|_{L_{2} ([0,1), (dt)/(1-t))}.
\end{eqnarray*}
Theorem~\ref{thmBesovSpacesNormEquiv} completes the proof, since in
the case that
$Y$ is the Brownian motion, we have $H_G(t,Y_t)=|D^2 F(t,W_t)|$ and
in the other case, we can use (\ref{eqnDerBMgBM}) and Theorem \ref
{thmBesovSpacesNormEquiv}
again to see that
\[
\bigl\|(1-t)^{(2-\theta)/2} \bigl\|
H_G(t,Y_t) \bigr\|_{L_{p}}
\bigr\|_{L_{2} ([0,1), (dt)/(1-t))} < \infty.
\]
\item[(b)] For all $0 < t \leq1$,
\begin{eqnarray*}
\bigl\| D^{Y,\theta}_1 g(Y_1)
\bigr\|_{L_p} &\ge& \biggl\|
\biggl( \int_0^t (1-s)^{1-\theta}
H_G^2 (s, Y_s) \,ds \biggr)
^{1/2} \biggr\|_{L_{p}}
\\
&\ge& (1-t)^{(1-\theta)/2} \biggl\|\biggl(
\int
_0^t H_G^2 (s,
Y_s) \,ds \biggr) ^{1/2} \biggr\|
_{L_{p}}.
\end{eqnarray*}
If $Y$ is the Brownian motion, then we can bound this from below by
\[
\frac{1}{c_{\fontsize{8.36pt}{10pt}\selectfont{(\ref{eqnBDG})}}} (1-t)^{(1-\theta)/2} \bigl\|\nabla
F(t,W_t)- \nabla F(0,W_0) \bigr\|
_{L_p},
\]
where we have used (\ref{eqnBDG}) and (\ref
{eqnrepresentationnablaF}). This implies that
\[
\bigl\| \nabla F(t,W_t) \bigr\|_{L_p} \le\bigl\| \nabla
F(0,W_0) \bigr\|_{L_p} + c_{\fontsize{8.36pt}{10pt}\selectfont{(\ref{eqnBDG})}} (1-t)^{(\theta-1)/2}
\bigl\| D^{W,\theta}_1 f(W_1) \bigr\|_{L_p}
\]
and Theorem~\ref{thmBesovSpacesNormEquiv} can be used again.
If $Y$ is the coordinate-wise geometric Brownian motion, then we get from
(\ref{eqnDerBMgBM}) that
\begin{eqnarray*}
& & \biggl\|\biggl( \int_0^t
H_G^2 (s, Y_s) \,ds \biggr)
^{1/2} \biggr\|_{L_{p}}
\\
&&\qquad \ge \biggl\|\biggl( \int_0^t
H_F^2 (s, W_s) \,ds \biggr)
^{1/2} \biggr\|_{L_{p}} - \biggl\|\biggl( \int_0^1 \bigl|\nabla
F(s, W_s)\bigr|^2 \,ds \biggr) ^{1/2} \biggr\|_{L_{p}}
\\
&&\qquad \ge \frac{1}{c_{\fontsize{8.36pt}{10pt}\selectfont{(\ref{eqnBDG})}}} \bigl\|\nabla
F(t,W_t)-
\nabla F(0,W_0) \bigr\|_{L_p} - \biggl\|\biggl( \int_0^1 \bigl|\nabla
F(s, W_s)\bigr|^2 \,ds \biggr) ^{1/2} \biggr\|_{L_{p}}
\\
&&\qquad \ge \frac{1}{c_{\fontsize{8.36pt}{10pt}\selectfont{(\ref{eqnBDG})}}} \bigl\|\nabla
F(t,W_t)
\bigr\|_{L_p} - \frac{1}{c_{\fontsize{8.36pt}{10pt}\selectfont{(\ref{eqnBDG})}}} \bigl\|\nabla F(0,W_0) \bigr\|_{L_p}
\\
&&\quad\qquad{}- \biggl\|\biggl( \int
_0^1 \bigl|\nabla F(s, W_s)\bigr|^2
\,ds \biggr) ^{1/2} \biggr\|_{L_{p}},
\end{eqnarray*}
where we again used the Burkholder--Davis--Gundy inequalities (\ref{eqnBDG}).
Because the last two terms on the right-hand side are finite,
we can conclude as in the case of the Brownian motion.

\item[(ii)] Because of (\ref{eqnDerBMgBM}) and
$ ( \int_0^1 |\nabla F(t,W_t)|^2 \,dt ) ^{1/2} \in
L_p$, we get
$D_1^{Y,1} g(Y_1)\in L_p$ if and only if
$ ( \int_0^1 | D^2F(t,W_t)|^2 \,dt ) ^{1/2} \in L_p$.
Using relations
(\ref{eqnd12hessian}) and (\ref{eqnrepresentationnablaF}),
one easily checks that this is equivalent to $f\in\mathbb{D}_{1,p}$.

\item[(iii)] Since (\ref{eqnDerBMgBM}) implies the equivalence of
\[
\bigl\| D^{Y,\theta}_1 g(Y_1)
\bigr\|_{L_2}^2 = \int
_0^1 (1-t)^{1-\theta} \bigl\|
H_G(t, Y_t) \bigr\|_{L_2}^2 \,dt <
\infty
\]
and $\int_0^1 (1-t)^{1-\theta} \| D^2F(t,W_t)\|_{L_2}^2 \,dt < \infty$,
we can use Theorem~\ref{thmBesovSpacesNormEquiv}.\quad\qed
\end{longlist}\noqed
\end{pf}


\section{An approximation problem in $L_p$}\label{secappoximation-p}

In the whole section, we use the convention (\ref{eqnconv1}) and
(\ref{eqnconv2}).

\subsection*{Time-nets}
Given a sequence $0=t_0 \le\cdots\le t_{n-1}<t_n=1$ and $0<\theta\le1$,
we let
\begin{eqnarray*}
\bigl|(t_i)_{i=0}^n\bigr|_\theta&:= & \sup
_{i=1,\ldots,n} \sup_{t_{i-1} \le u
< t_i} \frac{|t_{i}-u|}{(1-u)^{1-\theta}} = \sup
_{i=1,\ldots,n} \frac{|t_{i}-t_{i-1}|}{(1-t_{i-1})^{1-\theta}},
\\
\bigl|(t_i)_{i=0}^n\bigr| &:= & \bigl|(t_i)_{i=0}^n\bigr|_1
\end{eqnarray*}
so that $|(t_i)_{i=0}^n|$ is the usual mesh-size.
As special adapted deterministic time-nets we use
$\tau_n^\theta=(t_{i,n}^\theta)_{i=0}^n$ defined by
\[
t_{i,n}^\theta:= 1 - \biggl( 1 - \frac{i}{n} \biggr)
^{1/\theta}.
\]
For these time-nets,
%
\begin{eqnarray}\label{eqnestimatenets1}
\bigl|t_{i,n}^\theta-u\bigr| &\le&\frac{|t_{i,n}^\theta
-u|}{(1-u)^{1-\theta}} \le
\frac{|t_{i,n}^\theta- t_{i-1,n}^\theta|}{(1-t_{i-1,n}^\theta
)^{1-\theta}} \le\frac{1}{\theta n}
\nonumber\\[-10pt]\\[-10pt]
\eqntext{\mbox{for } u \in\bigl[t_{i-1,n}^\theta,t_{i,n}^\theta\bigr),}
\end{eqnarray}
which implies that
%
\begin{equation}
\label{eqnupperboundthetanormthetanet} \bigl|\tau_n^\theta\bigr| \le\bigl|\tau
_n^\theta\bigr|_\theta
\le\frac{1}{\theta n}.
\end{equation}
Moreover, we have that
%
\begin{eqnarray}
\label{eqnestimatenets2} \frac{(1-t_{i-1,n}^\theta)^{1-\theta
}}{|t_{i,n}^\theta-
t_{i-1,n}^\theta|} &\le& \beta n
\end{eqnarray}
for some $\beta>0$ independent from $n$.

\subsection*{The basic equivalence in $L_p$} The following result
reduces the computation of the
$L_p$-norm of the error processes defined in
Definition~\ref{deferrorprocesses}
to an expression involving the curvature $H_G(t,Y_t)$ similar to a
square function.
This result generalizes \cite{Gei15}, Theorem 4.4, proved for
deterministic nets in the
$L_2$-case.

%
\begin{theorem}\label{thmLp-error-equivalence}
For $2 \le p < \infty$, there\vspace*{2pt} is a constant
$c_{\fontsize{8.36pt}{10pt}\selectfont{(\ref{thmLp-error-equivalence} )}} \ge1$ depending
at most on $p$
such that for all $g(Y_1)\in L_p$ and
$\tau=(\tau_i)_{i=0}^n\in{\mathcal T}^{\mathrm{rand}}$ we have that
\begin{eqnarray*}
\bigl\| C_1\bigl(g(Y_1),\tau\bigr) \bigr\|_{L_p} &\le&
c_{\fontsize{8.36pt}{10pt}\selectfont{(\ref{thmLp-error-equivalence} )}} \Biggl\llVert\Biggl( \sum_{i=1}^n
\int_{\tau_{i-1}}^{\tau_i} (\tau_i-t)
H_G^2(t,Y_t) \,dt \Biggr) ^{1/2} \Biggr\rrVert_{L_p},
\\
\inf_v \bigl\| C_1\bigl(g(Y_1),\tau,v
\bigr) \bigr\|_{L_p} &\ge& \frac{1}{c_{\fontsize{8.36pt}{10pt}\selectfont{(\ref{thmLp-error-equivalence} )}}}
\Biggl\llVert\Biggl( \sum
_{i=1}^n \int_{\tau_{i-1}}^{\tau_i}
(\tau_i-t) H_G^2(t,Y_t) \,dt
\Biggr) ^{1/2} \Biggr\rrVert_{L_p},
\end{eqnarray*}
where the infimum is taken over all simple
random vectors $v_{\tau_{i-1}}\dvtx \Omega\to\mathbb{R}^d$ that are
$\mathcal{F}_{\tau_{i-1}}$-measurable.
\end{theorem}

%
\begin{remark}
\label{remarkthmLp-error-equivalence}
Both inequalities in Theorem~\ref{thmLp-error-equivalence} are
proved by stopping at $0<T<1$ and
letting $T\uparrow1$. Therefore, it might be possible for one or both
sides of an inequality to be
infinite. However, this cannot be the case: step (b) of our proof for
the trivial time-net $0=t_0<t_1=1$
gives by (\ref{eqnlowerboundapproximation}) that
\begin{eqnarray*}
& & \biggl\llVert\biggl( \int_0^T (T-t)
H_G^2(t,Y_t) \,dt \biggr) ^{1/2} \biggr\rrVert_{L_p}
\\
&&\qquad \le c \bigl\| \mathbb{E} \bigl(g(Y_1)|\mathcal{F}_T \bigr)
- \mathbb{E}g(Y_1) -\nabla G(0,Y_0)
(Y_T-Y_0) \bigr\|_{L_p}
\\
&&\qquad \le c \bigl\| g(Y_1) - \mathbb{E}g(Y_1) -\nabla
G(0,Y_0) (Y_1-Y_0) \bigr\|_{L_p} <
\infty
\end{eqnarray*}
so that
\begin{eqnarray*}
&& \Biggl\llVert\Biggl( \sum_{i=1}^n \int
_{\tau_{i-1}}^{\tau_i} (\tau_i-t)
H_G^2(t,Y_t) \,dt \Biggr) ^{1/2} \Biggr\rrVert_{L_p}
\\
&&\qquad \le\biggl\llVert\biggl( \int
_0^1 (1-t) H_G^2(t,Y_t)
\,dt \biggr) ^{1/2} \biggr\rrVert_{L_p} < \infty.
\end{eqnarray*}
Following (\ref{eqnupperboundapproximation}) from step (c), this
implies that
$\sup_{0\le T <1} \| C_T(g(Y_1),\tau) \|_{L_p}<\infty$, from which
we can conclude that
\[
\Biggl( \int_0^1 \sum
_{i=1}^n \chi_{(\tau_{i-1},\tau_i]}(t) \bigl|\nabla G(
\tau_{i-1},Y_{\tau_{i-1}}) \sigma(Y_t)\bigr|^2\,dt
\Biggr) ^{1/2} \in L_p
\]
and $C_1(g(Y_1),\tau)\in L_p$.
Finally, we have that
%
\begin{equation}
\label{eqnoptimalapproximationbetterthan-simple} \inf_v \bigl\| C_1
\bigl(g(Y_1),\tau,v\bigr) \bigr\|_{L_p}\le\bigl\| C_1
\bigl(g(Y_1),\tau\bigr) \bigr\|_{L_p},
\end{equation}
where the infimum is taken over all simple random vectors $v_{\tau
_{i-1}}\dvtx \Omega\to\mathbb{R}^d$
that are $\mathcal{F}_{\tau_{i-1}}$-measurable. The latter also
implies that all three expressions---in particular the simple and optimal $L_p$-approximation---in Theorem
\ref{thmLp-error-equivalence}
are equivalent up to a multiplicative constant.
\end{remark}

%
\begin{remark}
\label{remarklowerbound}
Theorem~\ref{thmLp-error-equivalence} provides an alternative way
to prove the lower bound
\[
\bigl\| C_1\bigl(g(Y_1),\tau^n\bigr)
\bigr\|_{L_p} \ge\frac{\delta}{\sqrt{n}}
\]
for some $\delta>0$, all $n=1,2,\ldots$ and all nets $\tau^n=(\tau
_i)_{i=1}^n\in{\mathcal T}^{\mathrm{rand}}$
whenever there is no row vector $v_0\in\mathbb{R}^d$ such that
$g(Y_1) = \mathbb{E}g(Y_1) + v_0 (Y_1-Y_0)$ a.s.
This lower bound was obtained in \cite{Gei20} in the one-dimensional case
using an asymptotic argument.
To check the lower bound, observe for
$0\le a < b \le1$ and $\rho_i:= (a\vee\tau_i) \wedge b$ that
\begin{eqnarray*}
&& \Biggl\llVert\Biggl( \sum_{i=1}^n \int
_{\tau_{i-1}}^{\tau_i} (\tau_i-t)
H_G^2(t,Y_t) \,dt \Biggr) ^{1/2} \Biggr\rrVert_{L_p}
\\
&&\qquad \ge \Biggl\llVert\Biggl( \sum
_{i=1}^n \int_{\rho_{i-1}}^{\rho_i}
(\rho_i-t) H_G^2(t,Y_t)
\,dt \Biggr) ^{1/2} \Biggr\rrVert_{L_p}
\\
&&\qquad \ge \Biggl\llVert\inf_{a\le t < b} H_G(t,Y_t)
\sum_{i=1}^n \frac{\rho_{i}-{\rho_{i-1}}}{\sqrt{2n}} \Biggr
\rrVert_{L_p}
\\
&&\qquad  =  \frac{b-a}{\sqrt{2n}} \Bigl\llVert\inf_{a\le t < b}
H_G(t,Y_t) \Bigr\rrVert_{L_p}.
\end{eqnarray*}
Assume now $\sup_{0\le a < b \le1} \llVert\inf_{a\le t < b} H_G(t,Y_t)
\rrVert_{L_p} = 0$ and fix $0<T<1$ and $0< a_n\uparrow T$. Then
\[
0 = \lim_n \Bigl\llVert\inf_{a_n \le t < T}
H_G(t,Y_t) \Bigr\rrVert_{L_p} = \Bigl\llVert
\lim_n \inf_{a_n \le t < T} H_G(t,Y_t)
\Bigr\rrVert_{L_p} = \bigl\llVert H_G(T,Y_T)
\bigr\rrVert_{L_p}
\]
so that
$H_G(T,Y_T) = 0$ a.s. for all $0<T<1$. Applying Theorem \ref
{thmLp-error-equivalence}
for the trivial time-net $\{ 0,1 \}$ yields
\[
\bigl\| g(Y_1) - \mathbb{E}g(Y_1) - \nabla
G(0,Y_0) (Y_1-Y_0) \bigr\|_{L_p} = 0.
\]
\end{remark}

\begin{pf*}{Proof of Theorem~\ref{thmLp-error-equivalence}}
(a) Assume a deterministic time $0<T<1$, two stopping times $0 \le a
\le b \le T$ and that $v_a$ is a
simple $\mathcal{F}_a$-measurable random (row) vector.
Exploiting relations (\ref{eqnrepresentationnablaF}) and (\ref
{eqnDerBMgBM})
one quickly checks that
\[
\biggl( \frac{\partial G}{\partial y_k}(b,Y_b) - v_a^k
\biggr) \sigma_{kk}(Y_b) = m_a(k) + \sum
_{l=1}^d \int_a^b
\lambda_u^a(k,l) \,dW_u^l
\]
with
\begin{eqnarray*}
m_a(k) &:= & \biggl( \frac{\partial G}{\partial y_k}(a,Y_a) -
v_a^k \biggr) \sigma_{kk}(Y_a)
\end{eqnarray*}
and
\begin{eqnarray*}
\lambda^a_u(k,l) &:= & \biggl( \sigma_{kk}
\sigma_{ll}\frac{\partial^2 G}{\partial
y_k\,\partial y_l} \biggr) (u,Y_u) + \biggl(
\frac{\partial G}{\partial y_l}(u,Y_u)-v_a^l \biggr)
\biggl( \sigma_{ll} \frac{\partial\sigma_{ll}}{\partial y_k} \biggr)
(Y_u),
\end{eqnarray*}
where $m_a:= (m_a(1),\ldots,m_a(d))$ will be considered as a row vector.

(b) \textit{Lower bound for} $\|C_1(g(Y_1),\tau,v)\|_{L_p}$:
Let us fix $0<T<1$ and define $\rho_i:=\tau_i\wedge T$ and $\alpha
_{\rho_{i-1}}:= v_{\tau_{i-1}}\chi_{\{\tau_{i-1}<T\}}$. Note that
$\rho_i$ and $\alpha_{\rho_{i-1}}$ are
$\mathcal{F}_{\rho_{i-1}}$-measurable for $i=1,\ldots,n$. Replacing
$v$ by $\alpha$ in the definitions of $m$ and $\lambda$ from
step (a), it follows that
\begin{eqnarray*}
C_T\bigl(g(Y_1),\tau,v\bigr) & = & \sum
_{i=1}^n \int_{\rho_{i-1}}^{\rho_i}
\biggl[ m_{\rho
_{i-1}} + \biggl[ \int_{\rho_{i-1}}^{t\vee\rho_{i-1}}
\lambda^{\rho_{i-1}}_u d W_u \biggr]^\top
\biggr] \,d W_t.
\end{eqnarray*}
Using (\ref{eqnBDG}) and the convexity inequality \cite{Garsia}, pp. 104--105,
p. 171, we achieve
\begin{eqnarray*}
\hspace*{-5pt}& & \Biggl\llVert\sum_{i=1}^n \int
_{\rho_{i-1}}^{\rho_i} \biggl[ m_{\rho_{i-1}} + \biggl[ \int
_{\rho_{i-1}}^{t\vee\rho_{i-1}} \lambda^{\rho_{i-1}}_u d
W_u \biggr]^\top\biggr] \,d W_t \Biggr\rrVert
_{L_p}^p
\\
\hspace*{-5pt}&&\qquad \ge c_{\fontsize{8.36pt}{10pt}\selectfont{(\ref{eqnBDG})}}^{-p} \mathbb{E} \Biggl( \sum
_{i=1}^n \int_{\rho_{i-1}}^{\rho_i}
\biggl\llvert m_{\rho_{i-1}} + \biggl[ \int_{\rho_{i-1}}^{t\vee\rho_{i-1}}
\lambda^{\rho_{i-1}}_u d W_u \biggr]^\top
\biggr\rrvert^2 \,d t \Biggr) ^{p/2}
\\
\hspace*{-5pt}&&\qquad \ge c_{\fontsize{8.36pt}{10pt}\selectfont{(\ref{eqnBDG})}}^{-p} (p/2)^{-p/2} \mathbb{E} \Biggl(
\sum_{i=1}^n \mathbb{E}_{\mathcal{F}_{\rho
_{i-1}}}
\int_{\rho_{i-1}}^{\rho_i} \biggl\llvert m_{\rho_{i-1}} +
\biggl[ \int_{\rho_{i-1}}^{t\vee\rho_{i-1}} \lambda^{\rho
_{i-1}}_u
d W_u \biggr]^\top\biggr\rrvert^2 \,d t
\Biggr) ^{p/2}
\\
\hspace*{-5pt}&&\qquad \ge (c_{\fontsize{8.36pt}{10pt}\selectfont{(\ref{eqnBDG})}} \sqrt{p/2})^{-p} \mathbb{E} \Biggl( \sum
_{i=1}^n \mathbb{E}_{\mathcal{F}_{\rho
_{i-1}}} \int
_{\rho_{i-1}}^{\rho_i} \llvert m_{\rho_{i-1}} \rrvert
^2 \,d t \Biggr) ^{p/2}
\\
\hspace*{-5pt}&&\qquad  =  (c_{\fontsize{8.36pt}{10pt}\selectfont{(\ref{eqnBDG})}} \sqrt{p/2})^{-p} \mathbb{E} \Biggl( \sum
_{i=1}^n \int_{\rho_{i-1}}^{\rho_i}
\llvert m_{\rho_{i-1}} \rrvert^2 \,d t \Biggr) ^{p/2},
\end{eqnarray*}
where we used the assumption that $\rho_i$ is $\mathcal{F}_{\rho
_{i-1}}$-measurable.
From this, we deduce that
\begin{eqnarray*}
& & \Biggl\llVert\sum_{i=1}^n \int
_{\rho_{i-1}}^{\rho_i} \biggl[ \int_{\rho_{i-1}}^{t\vee\rho_{i-1}}
\lambda^{\rho_{i-1}}_u d W_u \biggr]^\top
d W_t \Biggr\rrVert_{L_p}
\\
&&\qquad \le \Biggl\llVert\sum_{i=1}^n \int
_{\rho_{i-1}}^{\rho_i} m_{\rho
_{i-1}} \,d W_t \Biggr
\rrVert_{L_p}
\\
&&\quad\qquad{} + \Biggl\llVert\sum_{i=1}^n \int
_{\rho_{i-1}}^{\rho_i} \biggl[ m_{\rho_{i-1}} + \biggl[ \int
_{\rho_{i-1}}^{t\vee\rho_{i-1}} \lambda^{\rho_{i-1}}_u d
W_u \biggr]^\top\biggr] \,d W_t \Biggr\rrVert
_{L_p}
\\
&&\qquad \le c_{\fontsize{8.36pt}{10pt}\selectfont{(\ref{eqnBDG})}} \Biggl\llVert\Biggl( \sum_{i=1}^n
\int_{\rho
_{i-1}}^{\rho_i} |m_{\rho_{i-1}}|^2 \,d t
\Biggr) ^{1/2} \Biggr\rrVert_{L_p}
\\
&&\quad\qquad{} + \Biggl\llVert\sum_{i=1}^n \int
_{\rho_{i-1}}^{\rho_i} \biggl[ m_{\rho_{i-1}} + \biggl[ \int
_{\rho_{i-1}}^{t\vee\rho_{i-1}} \lambda^{\rho_{i-1}}_u d
W_u \biggr]^\top\biggr] \,d W_t \Biggr\rrVert
_{L_p}
\\
&&\qquad \le \bigl[c_{\fontsize{8.36pt}{10pt}\selectfont{(\ref{eqnBDG})}}^2 \sqrt{p/2} + 1\bigr] \Biggl\llVert
\sum_{i=1}^n \int_{\rho_{i-1}}^{\rho_i}
\biggl[ m_{\rho
_{i-1}} + \biggl[ \int_{\rho_{i-1}}^{t\vee\rho_{i-1}}
\lambda^{\rho_{i-1}}_u d W_u \biggr]^\top
\biggr] \,d W_t \Biggr\rrVert_{L_p}
\end{eqnarray*}
so that
\begin{eqnarray*}
&& \bigl\| C_T\bigl(g(Y_1),\tau,v\bigr) \bigr\|_{L_p}
\\
&&\qquad \ge
{\bigl[c_{\fontsize{8.36pt}{10pt}\selectfont{(\ref{eqnBDG})}}^2 \sqrt{p/2} + 1\bigr]}^{-1}
\Biggl\llVert\sum_{i=1}^n \int
_{\rho_{i-1}}^{\rho_i} \biggl[ \int_{\rho
_{i-1}}^{t\vee\rho_{i-1}}
\lambda^{\rho_{i-1}}_u d W_u \biggr]^\top
d W_t \Biggr\rrVert_{L_p}.
\end{eqnarray*}
We continue by writing
\begin{eqnarray*}
&& \Biggl\llVert\sum_{i=1}^n \int
_{\rho_{i-1}}^{\rho_i} \biggl[ \int_{\rho_{i-1}}^{t\vee\rho_{i-1}}
\lambda^{\rho_{i-1}}_u d W_u \biggr]^\top
d W_t \Biggr\rrVert_{L_p}
\\
&&\qquad  =  \Biggl\llVert\sum_{i=1}^n \int
_0^1 \biggl[ \int_0^t
\bigl[ \chi_{(\rho_{i-1},t\vee\rho_{i-1}]}(u) \chi_{(\rho_{i-1},\rho
_i]}(t) \lambda^{\rho_{i-1}}_u
\bigr] \,d W_u \biggr]^\top d W_t \Biggr\rrVert
_{L_p}
\\
&&\qquad  =  \biggl\llVert\int_0^1 \biggl[ \int
_0^t \mu^\rho(t,u) \,d
W_u \biggr]^\top d W_t \biggr\rrVert
_{L_p}
\end{eqnarray*}
with the $d\times d$-matrix
\[
\mu^\rho(t,u):= \sum_{i=1}^n
\chi_{(\rho_{i-1},t]} (u) \chi_{(\rho_{i-1},\rho
_i]}(t) \lambda^{\rho_{i-1}}_u
= \sum_{i=1}^n \chi_{\{ \rho_{i-1}<u\le t \le\rho_i \}}
\lambda^{\rho_{i-1}}_u.
\]
Here, we used again the condition that $\rho_i$ is $\mathcal{F}_{\rho
_{i-1}}$-measurable.
By (\ref{eqnBDG}) and Lemma~\ref{lemmadoubleIntBDG} (note that
$\rho_i\le T < 1$),
\begin{eqnarray*}
&& \biggl\llVert\int_0^1 \biggl[ \int
_0^t \mu^\rho(t,u) \,d
W_u \biggr]^\top d W_t \biggr\rrVert
_{L_p}
\\
&&\qquad \sim_{c_{\fontsize{6.6pt}{10pt}\selectfont{(\ref{eqnBDG})}}} \biggl\llVert\biggl
( \int
_0^1 \biggl\llvert\int_0^1
\mu^\rho(t,u) \,d W_u \biggr\rrvert^2 \,d t
\biggr) ^{1/2} \biggr\rrVert_{L_p}
\\
&&\qquad \sim_{c_{\fontsize{6.6pt}{10pt}\selectfont{(\ref{lemmadoubleIntBDG} )}}} \biggl\llVert
\biggl( \int
_0^1 \int_0^t
\bigl\llvert\mu^\rho(t,u) \bigr\rrvert^2 \,du \,dt
\biggr)^{1/2} \biggr\rrVert_{L_p}
\\
&&\qquad  =  \Biggl\llVert\Biggl( \sum_{i=1}^n
\int_{\rho
_{i-1}}^{\rho_i} (\rho_i-t) \bigl|
\lambda_t^{\rho_{i-1}}\bigr|^2 \,dt \Biggr) ^{1/2} \Biggr\rrVert_{L_p}.
\end{eqnarray*}
Letting $\delta=0$ if $Y=W$ and $\delta=1$ if $Y$ is the geometric
Brownian motion, this can be combined with
\begin{eqnarray*}
& & \Biggl\llVert\Biggl( \sum_{i=1}^n
\int_{\rho_{i-1}}^{\rho_i} (\rho_i-t)
H_G^2(t,Y_t) \,dt \Biggr) ^{1/2} \Biggr\rrVert_{L_p}
\\
&&\quad\qquad{} - \delta\Biggl\llVert\Biggl( \sum_{i=1}^n
\int_{\rho_{i-1}}^{\rho
_i} (\rho_i-t) \bigl\llvert
\bigl( \nabla G(t,Y_t) - \alpha_{\rho_{i-1}} \bigr)
\sigma(Y_t) \bigr\rrvert^2 \,dt \Biggr) ^{1/2} \Biggr\rrVert_{L_p}
\\
&&\qquad \le \Biggl\llVert\Biggl( \sum_{i=1}^n
\int_{\rho_{i-1}}^{\rho_i} (\rho_i-t) \bigl|
\lambda_t^{\rho_{i-1}}\bigr|^2 \,dt \Biggr) ^{1/2} \Biggr\rrVert_{L_p}
\\
&&\qquad \le \Biggl\llVert\Biggl( \sum_{i=1}^n
\int_{\rho_{i-1}}^{\rho_i} (\rho_i-t)
H_G^2(t,Y_t) \,dt \Biggr) ^{1/2} \Biggr\rrVert_{L_p}
\\
&&\quad\qquad{} + \delta\Biggl\llVert\Biggl( \sum_{i=1}^n
\int_{\rho_{i-1}}^{\rho
_i} (\rho_i-t) \bigl\llvert
\bigl( \nabla G(t,Y_t) - \alpha_{\rho_{i-1}} \bigr)
\sigma(Y_t) \bigr\rrvert^2 \,dt \Biggr) ^{1/2} \Biggr\rrVert_{L_p}
\end{eqnarray*}
so that
\begin{eqnarray*}
&& \bigl\| C_T\bigl(g(Y_1),\tau,v\bigr) \bigr\|_{L_p}
\\
&&\qquad \ge\bigl[c_{\fontsize{8.36pt}{10pt}\selectfont{(\ref{eqnBDG})}}^2 \sqrt{p/2} + 1\bigr]^{-1}
c_{\fontsize{8.36pt}{10pt}\selectfont{(\ref{eqnBDG})}}^{-1} c_{\fontsize{8.36pt}{10pt}\selectfont{(\ref{lemmadoubleIntBDG} )}}^{-1}
\\
&&\quad\qquad{}\times  \Biggl[ \Biggl
\llVert\Biggl( \sum_{i=1}^n \int
_{\rho_{i-1}}^{\rho_i} (\rho_i-t)
H_G^2(t,Y_t) \,dt \Biggr) ^{1/2} \Biggr\rrVert_{L_p}
\\
&&\hspace*{51pt}{} - \delta\Biggl\llVert\Biggl( \sum_{i=1}^n
\int_{\rho_{i-1}}^{\rho_i} \bigl\llvert\bigl( \nabla
G(t,Y_t) - \alpha_{\rho_{i-1}} \bigr) \sigma(Y_t) \bigr
\rrvert^2 \,dt \Biggr) ^{1/2} \Biggr\rrVert
_{L_p} \Biggr].
\end{eqnarray*}
In the case of the Brownian motion, the last term disappears. In the
case of the geometric Brownian motion,
we apply again (\ref{eqnBDG}) to see that
\[
\Biggl\llVert\Biggl( \sum_{i=1}^n \int
_{\rho_{i-1}}^{\rho_i} \bigl| \bigl(\nabla G(t,Y_t)-
\alpha_{\rho_{i-1}}\bigr) \sigma(Y_t)\bigr|^2 \,dt \Biggr)
^{1/2} \Biggr\rrVert_{L_p} \le c_{\fontsize{8.36pt}{10pt}\selectfont{(\ref{eqnBDG})}} \bigl\|
C_T\bigl(g(Y_1),\tau,v\bigr) \bigr\|_{L_p}.
\]
Hence, in both cases, we have that
%
\begin{eqnarray}\label{eqnlowerboundapproximation}
&& \inf_v \bigl\| C_T
\bigl(g(Y_1),\tau,v\bigr) \bigr\|_{L_p}\nonumber
\\
&&\qquad \ge\frac{1}{[c_{\fontsize{8.36pt}{10pt}\selectfont{(\ref{eqnBDG})}}^2 \sqrt{p/2} + 1] c_{\fontsize{8.36pt}{10pt}\selectfont{(\ref{eqnBDG})}} c_{\fontsize{8.36pt}{10pt}\selectfont{(\ref{lemmadoubleIntBDG} )}}+ c_{\fontsize{8.36pt}{10pt}\selectfont{(\ref{eqnBDG})} }}
\\
&&\quad\qquad{}\times  \Biggl\llVert\Biggl( \sum_{i=1}^n
\int_{\rho_{i-1}}^{\rho_i} (\rho_i-t)
H_G^2(t,Y_t) \,dt \Biggr) ^{1/2} \Biggr\rrVert_{L_p}.\nonumber
\end{eqnarray}
By $T\uparrow1$, we obtain the lower bound of our theorem.

(c) \textit{Upper bound for} $\|C_1(g(Y_1),\tau)\|_{L_p}$:
For $0<T<1$, using the arguments and notation from step (b) and
\[
\nu^{\rho_{i-1}}_u:= \bigl( \nabla G(u,Y_u)-
\nabla G(\tau_{i-1},Y_{\tau_{i-1}})\chi_{\{
\tau_{i-1}<T\}} \bigr)
\sigma(Y_u),
\]
we obtain
\begin{eqnarray*}
& & \bigl\| C_T\bigl(g(Y_1),\tau\bigr) \bigr\|_{L_p}
\\[-1pt]
&&\qquad \le c_{\fontsize{8.36pt}{10pt}\selectfont{(\ref{eqnBDG})}} c_{\fontsize{8.36pt}{10pt}\selectfont{(\ref{lemmadoubleIntBDG} )}} \Biggl
\llVert\Biggl( \sum
_{i=1}^n \int_{\rho_{i-1}}^{\rho_i}
(\rho_i-t) H_G^2(t,Y_t) \,dt
\Biggr) ^{1/2} \Biggr\rrVert_{L_p}
\\[-1pt]
&& \quad\qquad{} + c_{\fontsize{8.36pt}{10pt}\selectfont{(\ref{eqnBDG})}} c_{\fontsize{8.36pt}{10pt}\selectfont{(\ref{lemmadoubleIntBDG}
)}}\delta\Biggl\llVert\Biggl( \int
_0^T \int_0^t
\Biggl\llvert\sum_{i=1}^n
\chi_{\{\rho_{i-1}<u\le t \le\rho_i\}} \nu^{\rho_{i-1}}_u \Biggr\rrvert
^2 \,d u\, d t \Biggr) ^{1/2} \Biggr\rrVert
_{L_p}
\\[-1pt]
&&\qquad \le c_{\fontsize{8.36pt}{10pt}\selectfont{(\ref{eqnBDG})}} c_{\fontsize{8.36pt}{10pt}\selectfont{(\ref{lemmadoubleIntBDG} )}} \Biggl
\llVert\Biggl( \sum
_{i=1}^n \int_{\tau_{i-1}}^{\tau_i}
(\tau_i-t) H_G^2(t,Y_t) \,dt
\Biggr) ^{1/2} \Biggr\rrVert_{L_p}
\\[-1pt]
&&\quad\qquad{} + c_{\fontsize{8.36pt}{10pt}\selectfont{(\ref{eqnBDG})}} c_{\fontsize{8.36pt}{10pt}\selectfont{(\ref{lemmadoubleIntBDG}
)}}\delta\Biggl\llVert\Biggl( \int
_0^T \int_0^t
\sum_{i=1}^n \chi_{(\rho
_{i-1},\rho_i]}(u) \bigl|
\nu^{\rho_{i-1}}_u \bigr|^2 \,d u\, d t \Biggr)
^{1/2} \Biggr\rrVert_{L_p}.
\end{eqnarray*}
Because $2\le p < \infty$, we can continue by
\begin{eqnarray*}
&& \Biggl\| \Biggl( \int_0^T \int
_0^t \sum_{i=1}^n
\chi_{(\rho_{i-1},\rho_i]}(u) \bigl\llvert\nu^{\rho_{i-1}}_u \bigr
\rrvert^2 \,d u\, d t \Biggr)^{1/2} \Biggr\|_{L_p}
\\[-1pt]
&&\qquad \le \Biggl( \int_0^T \Biggl\| \Biggl( \int
_0^t \sum_{i=1}^n
\chi_{(\rho_{i-1},\rho_i]}(u) \bigl\llvert\nu^{\rho_{i-1}}_u \bigr
\rrvert^2 \,d u \Biggr)^{1/2} \Biggr\|_{L_p}^2
\,d t \Biggr)^{1/2}
\\[-1pt]
&&\qquad \le c_{\fontsize{8.36pt}{10pt}\selectfont{(\ref{eqnBDG})}} \Biggl( \int_0^T \Biggl\|
\int_0^t \sum_{i=1}^n
\chi_{(\rho_{i-1},\rho_i]}(u) \nu^{\rho_{i-1}}_u d W_u
\Biggr\|_{L_p}^2 \,d t \Biggr)^{1/2}
\\[-1pt]
&&\qquad  =  c_{\fontsize{8.36pt}{10pt}\selectfont{(\ref{eqnBDG})}} \biggl( \int_0^T \bigl\|
C_t\bigl(g(Y_1),\tau\bigr) \bigr\|_{L_p}^2
\,d t \biggr)^{1/2}.
\end{eqnarray*}
Combining these estimates, we achieve
\begin{eqnarray*}
&& \bigl\| C_T\bigl(g(Y_1),\tau\bigr) \bigr\|_{L_p}^2
\\[-2pt]
&&\qquad \le2 c_{\fontsize{8.36pt}{10pt}\selectfont{(\ref{eqnBDG})}}^2 c_{\fontsize{8.36pt}{10pt}\selectfont{(\ref{lemmadoubleIntBDG} )}}^2 \Biggl
\llVert\Biggl( \sum_{i=1}^n \int
_{\tau_{i-1}}^{\tau_i} (\tau_i-t)
H_G^2(t,Y_t) \,dt \Biggr) ^{1/2} \Biggr\rrVert_{L_p}^2
\\[-1pt]
&&\quad\qquad{}
+ 2 c_{\fontsize{8.36pt}{10pt}\selectfont{(\ref{eqnBDG})}}^4 c_{\fontsize{8.36pt}{10pt}\selectfont{(\ref{lemmadoubleIntBDG}
)}}^2 \int
_0^T \bigl\| C_t\bigl(g(Y_1),
\tau\bigr) \bigr\|_{L_p}^2 \,dt.
\end{eqnarray*}
Gronwall's lemma thus implies that
%
\begin{eqnarray}\label{eqnupperboundapproximation}
&& \bigl\| C_T\bigl(g(Y_1),\tau\bigr)
\bigr\|_{L_p}
\nonumber\\[-3pt]\\[-19pt]
&&\qquad \le\sqrt{2} c_{\fontsize{8.36pt}{10pt}\selectfont{(\ref{eqnBDG})}} c_{\fontsize{8.36pt}{10pt}\selectfont{(\ref{lemmadoubleIntBDG} )}} e^{c_{\fontsize{8.36pt}{10pt}\selectfont{(\ref{eqnBDG})}}^4 c_{\fontsize{8.36pt}{10pt}\selectfont{(\ref{lemmadoubleIntBDG} )}}^2} \Biggl\llVert
\Biggl( \sum_{i=1}^n \int
_{\tau_{i-1}}^{\tau_i} (\tau_i-t)
H_G^2(t,Y_t) \,dt \Biggr) ^{1/2} \Biggr\rrVert_{L_p}.\nonumber
\end{eqnarray}
Finally, by $T\uparrow1$ we obtain the upper bound in Theorem \ref
{thmLp-error-equivalence}.
\end{pf*}
%

%
\begin{remark}\label{remarkrestrictiontimenets}
Our proof requires the assumption that the stopping time $\tau_i$ is
$\mathcal{F}_{\tau_{i-1}}$-measurable
so that $\rho_i$ is $\mathcal{F}_{\rho_{i-1}}$-measurable.
For example, we need that the field $(\mu^\rho(t,u))_{t,u\in[0,1]}$
has the property that
$\mu^\rho(t,u)$ is $\mathcal{F}_u$-measurable.
Moreover, in step (b) we used
$ \mathbb{E}_{\mathcal{F}_{\rho_{i-1}}} \int_{\rho_{i-1}}^{\rho
_i} \llvert m_{\rho_{i-1}} \rrvert^2 \,dt
= \break \int_{\rho_{i-1}}^{\rho_i} \llvert m_{\rho_{i-1}} \rrvert^2 \,dt$.
\end{remark}

\subsection*{Approximation with adapted time-nets in $L_p$}
We recall that the nets $\tau_n^\theta$ are given by
\[
t_{i,n}^\theta= 1 - \biggl( 1-\frac{i}{n} \biggr)
^{1/\theta}.
\]
The following result extends \cite{Gei19}, Theorem 3.2, from the
one-dimensional \mbox{$L_2$-}setting,
but see also \cite{Hujo2007}, Theorem 1, for a related $d$-dimensional
$L_2$-result.

%
\begin{theorem}\label{thmLp-adaptedtimenets}
For $2 \le p < \infty$, $0<\theta\le1$ and $g(Y_1)\in L_p$, the
following assertions are equivalent:
\begin{longlist}[(iii)]
\item[(i)]$\| D_1^{Y,\theta}g(Y_1)\|_{L_p} < \infty$.\vspace*{3pt}
\item[(ii)]$\sup_{\tau\in{\mathcal T}^{\mathrm{rand}}}
\frac{\| C_1(g(Y_1),\tau) \|_{L_p}}{\| \sqrt{|\tau|_\theta} \|
_{L_\infty}} < \infty$.\vspace*{3pt}
\item[(iii)]$\sup_{n\ge1} \sqrt{n} \| C_1(g(Y_1),\tau_n^\theta) \|_{L_p}
< \infty$.
\end{longlist}
In particular, for all $\tau\in{\mathcal T}^{\mathrm{rand}}$,
%
\begin{equation}
\label{eqnthmLp-adaptedtimenets} \bigl\| C_1\bigl(g(Y_1),\tau\bigr)
\bigr\|_{L_p} \le c_{\fontsize{8.36pt}{10pt}\selectfont{(\ref{thmLp-error-equivalence})}}\bigl\| \sqrt{|\tau
|_\theta} D_1^{Y,\theta}g(Y_1)\bigr\|_{L_p},
\end{equation}
where $c_{\fontsize{8.36pt}{10pt}\selectfont{(\ref{thmLp-error-equivalence})}}\ge1$ is the constant from
Theorem~\ref{thmLp-error-equivalence}.
\end{theorem}

For the proof, we need the following lemma that extends \cite{Gei19}, Lemma~3.8.

%
\begin{lemma}\label{lemmanon-uniformtime-nets}
Let $0< \theta\leq1$ and $0< p < \infty$.
Assume that $(\phi_t)_{t\in[0,1)}$ is a measurable process where all
paths are continuous and nonnegative.
Then the following assertions are equivalent:
\begin{longlist}[(iii)]
\item[(i)]
There exists a constant $c_1 >0$ such that
\[
\Biggl\|\sum_{i=1}^n
\int_{t_{i-1}}^{t_i} (t_i -u)
\phi_u \,du \Biggr\|_{L_{p}} \leq
c_1 \sup_{1 \leq i \leq n} \frac{t_i -
t_{i-1}}{(1-t_{i-1})^{1-\theta}}
\]
for all deterministic time-nets $0= t_0 <t_1 < \cdots< t_n = 1$.
\item[(ii)]
There exists a constant $c_2 >0$ such that, for all $n=1,2,\ldots,$
\[
\Biggl\|\sum_{i=1}^n
\int_{t_{i-1,n}^{\theta
}}^{t_{i,n}^{\theta}} \bigl(t^{\theta}_{i,n}
-u\bigr) \phi_u \,du \Biggr\|_{L_{p}}
\leq\frac{c_2}{n}.
\]
\item[(iii)]
There exists a constant $c_3 >0$ such that
\[
\biggl\|\int_0^1
(1-u)^{1-\theta} \phi_u \,du \biggr\|
_{L_{p}} \leq c_3.
\]
\end{longlist}
\end{lemma}

\begin{pf}
The implications (iii)${}\Rightarrow{}$(i)${}\Rightarrow{}$(ii)
are similar to \cite{Gei19}, Lemma 3.8. For (ii)${}\Rightarrow{}$(iii), take a sequence of deterministic nets $\tau^n =
(t_i^n)_{i=0}^n $ with
$0=t_0^n < t_1^n < \cdots< t_n^n = 1$ such that
\[
\bigl|\tau^n\bigr| \le\frac{\alpha}{n} \quad\mbox{and}\quad\sup
_{1\le i \le n} \frac{ (1-t_{i-1}^n)^{1-\theta}}{t_i^n -
t_{i-1}^n} \le\beta n
\]
for some $\alpha,\beta>0$ independent from $n$
[see, e.g., (\ref{eqnestimatenets1}) and (\ref
{eqnestimatenets2})]. For a fixed $0<T<1$, we define
\[
N_T^n:= \bigl\{ i \in\{ 1,\ldots,n \}\dvtx
t_{i-1}^n < T \bigr\}
\]
and observe that
\begin{eqnarray*}
&& \int_0^T (1-u)^{1-\theta}
\phi_u \,du
\\
&&\qquad \le\liminf_{n \to\infty} \sum
_{i\in N_T^n} \bigl(1-t_{i-1}^n
\bigr)^{1-\theta
} \phi_{t_{i-1}^n} \bigl(t_i^n-
t_{i-1}^n\bigr)
\end{eqnarray*}
for all $\omega\in\Omega$ because $\phi$ is continuous on $[0,T]$.
Hence,
\begin{eqnarray*}
& & \biggl\|\int_0^T
(1-u)^{1-\theta} \phi_u \,du \biggr\|
_{L_{p}}
\\
&&\qquad \le \biggl\| \liminf_{n\to\infty} \biggl[ \sup_{1 \leq i \leq n}
\frac{(1-t_{i-1}^n)^{1-\theta}}{t_i^n
- t_{i-1}^n} \biggr] \biggl[ \sum_{i \in N_T^n}
\bigl(t_i^n - t_{i-1}^n
\bigr)^2 \phi_{t_{i-1}^n} \biggr] \biggr\|_{L_{p}}
\\
&&\qquad \le \beta\biggl\|\liminf_{n\to\infty} n
\biggl[ \sum_{i\in N_T^n} \bigl(t_i^n
- t_{i-1}^n\bigr)^2 \phi_{t_{i-1}^n}
\biggr]\biggr\|_{L_{p}}.
\end{eqnarray*}
Noticing that $(t_i^n - t_{i-1}^n)^2 = 2 \int_{t_{i-1}^n}^{t_i^n}
(t_i^n - u) \,du$ we continue with
\begin{eqnarray*}
& & \beta\biggl\|\liminf_{n\to\infty} n \biggl[ 2 \sum
_{i\in N_T^n} \int_{t_{i-1}^n}^{t_i^n}
\bigl(t_i^n - u\bigr) \,du\,\phi_{t_{i-1}^n} \biggr]
\biggr\|_{L_{p}}
\\
&&\qquad \le \beta\biggl\|\liminf_{n\to\infty} n \biggl[ 2 \sum
_{i\in N_T^n} \int_{t_{i-1}^n}^{t_i^n}
\bigl(t_i^n - u\bigr) \phi_u \,du
\\
&&\hspace*{85pt}{}+ \sum_{i\in N_T^n} \sup
_{t_{i-1}^n \leq u < t_i^n} \bigl|\phi_u - \phi_{t_{i-1}^n}\bigr|
\bigl(t_i^n - t_{i-1}^n
\bigr)^2 \biggr] \biggr\|_{L_{p}}
\\
&&\qquad \le \beta\biggl\|\liminf_{n\to\infty} \biggl[ 2 n \sum
_{i\in N_T^n} \int_{t_{i-1}^n}^{t_i^n}
\bigl(t_i^n - u\bigr) \phi_u \,du
\\
&&\hspace*{77pt}{} + \alpha\sup_{i\in N_T^n} \sup
_{t_{i-1}^n \leq u < t_i^n} \bigl|\phi_u - \phi_{t_{i-1}^n}\bigr| \biggr]
\biggr\|_{L_{p}}
\\
&&\qquad \le \beta\biggl\| 2 \liminf_{n\to\infty} n \sum
_{i\in N_T^n} \int_{t_{i-1}^n}^{t_i^n}
\bigl(t_i^n - u\bigr) \phi_u \,du
\\
&&\hspace*{43pt}{}
+ \alpha\limsup_{n\to\infty} \sup
_{i\in N_T^n} \sup_{t_{i-1}^n
\leq u < t_i^n} \bigl|\phi_u -
\phi_{t_{i-1}^n}\bigr| \biggr\|_{L_{p}}
\\
&&\qquad =  2 \beta\biggl\| \liminf_{n\to\infty} n \sum
_{i\in N_T^n} \int_{t_{i-1}^n}^{t_i^n}
\bigl(t_i^n - u\bigr) \phi_u \,du
\biggr\|_{L_{p}}
\\
&&\qquad \le 2 \beta\liminf_{n\to\infty} n \biggl\| \sum
_{i=1}^n \int_{t_{i-1}^n}^{t_i^n}
\bigl(t_i^n - u\bigr) \phi_u \,du \biggr\|
_{L_{p}},
\end{eqnarray*}
where we used Fatou's lemma.
Finally, by monotone convergence this implies that
\[
\biggl\|\int_0^1
(1-u)^{1-\theta} \phi_u \,du \biggr\|
_{L_{p}} \le2 \beta\liminf_{n\to\infty} n \Biggl\| \sum
_{i=1}^n \int_{t_{i-1}^n}^{t_i^n}
\bigl(t_i^n - u\bigr) \phi_u \,du
\Biggr\|_{L_{p}}.
\]\upqed
\end{pf}

\begin{pf*}{Proof of Theorem~\ref{thmLp-adaptedtimenets}}
First, we employ Theorem~\ref{thmLp-error-equivalence} to confirm
equation~(\ref{eqnthmLp-adaptedtimenets}) by
\begin{eqnarray*}
&& \bigl\| C_1\bigl(g(Y_1),\tau\bigr) \bigr\|_{L_p}
\\
&&\qquad \le c_{\fontsize{8.36pt}{10pt}\selectfont{(\ref{thmLp-error-equivalence})}} \Biggl\llVert\Biggl( \sum_{i=1}^n
\int_{\tau_{i-1}}^{\tau_i} (\tau_i -t)
H_G^2(t,Y_t) \,dt \Biggr) ^{1/2} \Biggr\rrVert_{L_p}
\\
&&\qquad \le c_{\fontsize{8.36pt}{10pt}\selectfont{(\ref{thmLp-error-equivalence})}} \Biggl\llVert\sqrt{|\tau
|_\theta} \Biggl( \sum
_{i=1}^n \int_{\tau
_{i-1}}^{\tau_i}
(1-t)^{1-\theta} H_G^2(t,Y_t) \,dt
\Biggr) ^{1/2} \Biggr\rrVert_{L_p}.
\end{eqnarray*}
Part (i)${}\Rightarrow{}$(ii) follows from (\ref
{eqnthmLp-adaptedtimenets}) and
part (ii)${}\Rightarrow{}$(iii) from $|\tau_n^\theta|_\theta\le
\frac{1}{\theta n}$
[see (\ref{eqnupperboundthetanormthetanet})].
To show that (iii)${}\Rightarrow{}$(i), we apply
Theorem~\ref{thmLp-error-equivalence} and (\ref
{eqnoptimalapproximationbetterthan-simple}) to see that
\[
\frac{c}{\sqrt{n}} \ge\bigl\| C_1\bigl(g(Y_1),
\tau_n^\theta\bigr)\bigr\|_{L_p} \ge\frac{1}{c_{\fontsize{8.36pt}{10pt}\selectfont{(\ref{thmLp-error-equivalence})}}}
\Biggl\llVert\Biggl( \sum_{i=1}^n \int
_{t_{i-1,n}^\theta
}^{t_{i,n}^\theta} \bigl(t_{i,n}^\theta-t
\bigr) H_G^2(t,Y_t) \,dt \Biggr)
^{1/2} \Biggr\rrVert_{L_p}.
\]
Lemma~\ref{lemmanon-uniformtime-nets} completes the proof.
\end{pf*}

\subsection*{Approximation with equidistant time-nets in $L_p$}
Here, we extend the \mbox{$L_2$-}re\-sults \cite{Gei18}, Theorem 2.3, and \cite{Gei19}, Theorem 3.5 for $q=\infty$,
to the $L_p$-case, as well as \cite{Toivola1}, Theorem 1.2, which
concerned deterministic time-nets and
the one-dimensional Brownian motion, to random time-nets and the
geometric Brownian motion.

%
\begin{theorem}\label{thmLp-non-adaptedtimenets}
For $2 \le p < \infty$, $0<\theta<1$ and $g(Y_1)\in L_p$, the
following assertions are equivalent:
\begin{longlist}[(iii)]
\item[(i)]$f \in\mathbb{B}^{\theta}_{p,\infty}$.\vspace*{3pt}
\item[(ii)]$\sup_{\tau\in{\mathcal T}^{\mathrm{rand}}} \frac{\|
C_1(g(Y_1);\tau) \|_{L_p}}{\||\tau|^{\theta/2} \|_{L_\infty}} <\infty$.\vspace*{3pt}
\item[(iii)]$\sup_{n=1,2,\ldots} n^{\theta/2} \| C_1(g(Y_1);\tau
_n) \|_{L_p}<\infty$, where
$\tau_n=(i/n)_{i=0}^n$ are the equidistant time-nets.
\end{longlist}
In particular, for $\frac{1}{p}=\frac{1}{q}+\frac{1}{r}$ with $p\le
q,r \le\infty$
and for all $\tau\in{\mathcal T}^{\mathrm{rand}}$,
%
\begin{eqnarray}\label{eqnthmLp-non-adaptedtimenets}
\bigl\| C_1\bigl(g(Y_1),\tau\bigr)
\bigr\|_{L_p}
&\le& c_{\fontsize{8.36pt}{10pt}\selectfont{(\ref{thmLp-error-equivalence} )}} \biggl( \int_0^1 \bigl\|
\sqrt{\psi(t)} \bigr\|_{L_q}^2 (1-t)^{\theta-2} \,dt \biggr)
^{1/2}
\nonumber\\[-9pt]\\[-9pt]
&&{} \times \sup_{t\in[0,1)} (1-t)^{1-(\theta/2)} \bigl\|
H_G(t,Y_t) \bigr\| _{L_r},\nonumber
\end{eqnarray}
where $c_{\fontsize{8.36pt}{10pt}\selectfont{(\ref{thmLp-error-equivalence} )}} \ge1$ is
the constant from
Theorem~\ref{thmLp-error-equivalence} and
\[
\psi(t,\omega):= \Bigl( \max_{i=1,\ldots,n} \bigl|\tau_i(
\omega)-\tau_{i-1}(\omega)\bigr| \Bigr) \wedge(1-t).
\]
\end{theorem}
%

%
\begin{remark}
The order for the equidistant nets can also be obtained
from Theorem~\ref{thmLp-adaptedtimenets} under the condition
$\| D_1^{Y,\theta}g(Y_1)\|_{L_p} < \infty$ because\break 
$|(i/n)_{i=0}^n|_\theta= n^{-\theta}$.
\end{remark}

\begin{pf*}{Proof of Theorem~\ref{thmLp-non-adaptedtimenets}}
To verify (\ref{eqnthmLp-non-adaptedtimenets}), we use Theorem
\ref{thmLp-error-equivalence}
and derive that
\begin{eqnarray*}
\bigl\| C_1\bigl(g(Y_1),\tau\bigr) \bigr\|_{L_p} &\le&
c_{\fontsize{8.36pt}{10pt}\selectfont{(\ref{thmLp-error-equivalence})}} \Biggl\llVert\Biggl( \sum_{i=1}^n
\int_{\tau_{i-1}}^{\tau_i} (\tau_i -t)
H_G^2(t,Y_t) \,dt \Biggr) ^{1/2} \Biggr\rrVert_{L_p}
\\[-1pt]
&\le& c_{\fontsize{8.36pt}{10pt}\selectfont{(\ref{thmLp-error-equivalence})}} \biggl\llVert\biggl( \int_0^1
\psi(t) H_G^2(t,Y_t) \,dt \biggr)
^{1/2} \biggr\rrVert_{L_p}
\\[-1pt]
&\le& c_{\fontsize{8.36pt}{10pt}\selectfont{(\ref{thmLp-error-equivalence})}} \biggl( \int_0^1 \bigl\|
\sqrt{\psi(t)} H_G(t,Y_t)\bigr\|_{L_p}^2
\,dt \biggr) ^{1/2}
\\[-1pt]
&\le& c_{\fontsize{8.36pt}{10pt}\selectfont{(\ref{thmLp-error-equivalence})}} \biggl( \int_0^1 \bigl\|
\sqrt{\psi(t)} \bigr\|_{L_q}^2 \bigl\| H_G(t,Y_t)
\bigr\| _{L_r}^2 \,dt \biggr) ^{1/2}
\\[-1pt]
&\le& c_{\fontsize{8.36pt}{10pt}\selectfont{(\ref{thmLp-error-equivalence})}} \biggl( \int_0^1 \bigl\|
\sqrt{\psi(t)} \bigr\|_{L_q}^2 (1-t)^{\theta-2} \,dt \biggr)
^{1/2}
\\[-1pt]
& &{}\times \sup_{t\in[0,1)} (1-t)^{1-(\theta/2)} \bigl\|
H_G(t,Y_t) \bigr\|_{L_r}.
\end{eqnarray*}
Part (i)${}\Rightarrow{}$(ii): we first observe that for
\[
|\tau|(\omega) = \max_{i=1,\ldots,n} \bigl|\tau_i(\omega)-\tau
_{i-1}(\omega)\bigr|
\]
we can compute (for $q=\infty$)
\begin{eqnarray*}
& & \int_0^1 \bigl\| \sqrt{\psi(t)}
\bigr\|_{L_\infty}^2 (1-t)^{\theta-2} \,dt
\\
&&\qquad =  \||\tau|\|_{L_\infty} \int_0^{1-\||\tau|\|_{L_\infty}}
(1-t)^{\theta-2} \,dt + \int_{1-\||\tau|\|_{L_\infty}}^1
(1-t)^{\theta-1} \,dt
\\
&&\qquad  =  \||\tau|\|_{L_\infty} \frac{1}{1-\theta} \bigl(\||\tau|\|
_{L_\infty}^{\theta-1} - 1\bigr) + \frac{1}{\theta} \||\tau|\|_{L_\infty}^\theta
\\
&&\qquad \le \frac{1}{\theta(1-\theta)} \||\tau|\|_{L_\infty}^\theta,
\end{eqnarray*}
so that letting
$q=\infty$ and $r=p$ in (\ref{eqnthmLp-non-adaptedtimenets}) we obtain
\[
\bigl\| C_1\bigl(g(Y_1),\tau\bigr) \bigr\|_{L_p} \le
\frac{c_{\fontsize{8.36pt}{10pt}\selectfont{(\ref{thmLp-error-equivalence} )}}}{\sqrt{\theta(1-\theta
)}} \||\tau|\|_{L_\infty}^{\theta/2}
\sup_{t\in[0,1)} (1-t)^{1-(\theta/2)} \bigl\| H_G(t,Y_t)
\bigr\| _{L_p}.
\]
It remains to check that
\[
\sup_{t\in[0,1)} (1-t)^{1-(\theta/2)} \bigl\| H_G(t,Y_t)
\bigr\|_{L_p} < \infty,
\]
whenever $f \in\mathbb{B}^{\theta}_{p,\infty}$.
This follows from Theorem~\ref{thmBesovSpacesNormEquiv}, where we
additionally use~(\ref{eqnDerBMgBM}) and the a priori estimate
%
\begin{equation}
\label{eqna-priori-bound-nabla} \sup_{t\in[0,1)} (1-t)^{1/2} \bigl\|
\nabla F (t,W_t) \bigr\|_{L_p} < \infty
\end{equation}
from Lemma~\ref{lemmaestimateDerFbyf} if $Y$ is the geometric
Brownian motion.

The implication (ii)${}\Rightarrow{}$(iii) is trivial.

Part (iii)${}\Rightarrow{}$(i): employing Theorem \ref
{thmLp-error-equivalence} and (\ref
{eqnoptimalapproximationbetterthan-simple}), we achieve
\begin{eqnarray*}
c n^{-\theta/2} &\ge& \bigl\| C_1\bigl(g(Y_1),
\tau_n\bigr) \bigr\|_{L_p}
\\
&\ge& \frac{1}{c_{\fontsize{8.36pt}{10pt}\selectfont{(\ref{thmLp-error-equivalence})}}} \Biggl\llVert
\Biggl( \sum
_{i=1}^n \int_{(i-1)/n}^{i/n}
\biggl( \frac{i}{n}-t \biggr) H_G^2(t,Y_t)
\,dt \Biggr) ^{1/2} \Biggr\rrVert_{L_p}
\\
&\ge& \frac{1}{c_{\fontsize{8.36pt}{10pt}\selectfont{(\ref{thmLp-error-equivalence})}}} \biggl\llVert
\biggl( \int_{(n-1)/n}^1
( 1-t ) H_G^2(t,Y_t) \,dt \biggr)
^{1/2} \biggr\rrVert_{L_p}
\\
&\ge& \frac{1}{c_{\fontsize{8.36pt}{10pt}\selectfont{(\ref{thmLp-error-equivalence})}}} \biggl\llVert
\biggl( \int_{(n-1)/n}^1
( 1-t ) H_G^2 \biggl( 1-\frac{1}{n},Y_{1-(1/n)}
\biggr) \,dt \biggr) ^{1/2} \biggr\rrVert_{L_p}
\\
& = & \frac{1}{c_{\fontsize{8.36pt}{10pt}\selectfont{(\ref{thmLp-error-equivalence})}}} \sqrt{\frac
{1}{2}} \frac{1}{n} \biggl
\llVert H_G \biggl( 1-\frac{1}{n},Y_{1-(1/n)} \biggr)
\biggr\rrVert_{L_p},
\end{eqnarray*}
where we use in the last inequality the martingale property of the processes
\[
\biggl( \biggl( \sigma_{kk} \sigma_{ll}\frac{\partial^2 G}{\partial
y_k\,\partial y_l}
\biggr) (t,Y_t) \biggr) _{t\in[0,1)}.
\]
The estimate above means that
\[
\biggl\llVert H_G \biggl( 1-\frac{1}{n},Y_{1-(1/n)}
\biggr) \biggr\rrVert_{L_p} \le\sqrt{2} c c_{\fontsize{8.36pt}{10pt}\selectfont{(\ref{thmLp-error-equivalence})}}
n^{1-(\theta/2)}
\]
for all $n = 2,3, \ldots.$
Consequently,
\[
\bigl\llVert H_G ( t,Y_t ) \bigr\rrVert
_{L_p} \le2^{1-(\theta/2)} \sqrt{2} c c_{\fontsize{8.36pt}{10pt}\selectfont{(\ref{thmLp-error-equivalence})}}
(1-t)^{(\theta/2)-1},
\]
which follows from the monotonicity of $\| H_G(t,Y_t) \|_{L_p}$.
Theorem~\ref{thmBesovSpacesNormEquiv} completes the proof, where we
use (\ref{eqna-priori-bound-nabla})
again.
\end{pf*}


\section{Further extensions}\label{secextensions}

We see different open questions and possible extensions, and briefly indicate
some of them here: first, one should clarify whether Theorem \ref
{thmLp-error-equivalence} holds
true without the additional assumption on the stopping times that $\tau
_i$ is
$\mathcal{F}_{\tau_{i-1}}$-measurable.
Second, the investigation to what extend the results of this paper can
be extended to path dependent terminal conditions $g(Y_{r_1},\ldots,Y_{r_L})$ and their limits would possibly require new techniques and
yield to a deeper insight into the approximation problem (cf. \cite{Geiss-Geiss-Gobet2012}).
Finally, an extension to more general diffusions would be of interest,
but might require a modification of the Besov spaces (see \cite{Gei18}) and a comparison
of these modified spaces to the spaces we have used in this paper.
As described in Remark~\ref{remarkOU-semigroup}, Proposition~\ref{propgeneralinterpol} below,
which does not relay on semi-groups, might be useful in this respect.


%
\begin{appendix}
\section*{Appendix}

A key step in the proof of Theorem~\ref{thmLp-error-equivalence} is
the following-known formulation of the Burkholder--Davis--Gundy inequalities.

%
\begin{lemma}\label{lemmadoubleIntBDG}
Assume that $\mu\dvtx  [0,1]\times[0,1] \times\Omega\to\mathbb{R}^{d
\times d}$ satisfies
the following assumptions:
\begin{longlist}[(iii)]
\item[(i)]$\mu\dvtx  [0,1]\times[0,u] \times\Omega\to\mathbb{R}^{d \times
d}$ is
$\mathcal{B}([0,1])\times\mathcal{B}([0,u])\times\mathcal
{F}_u$-measurable for\vspace*{1pt} all $u\in[0,1]$.
\item[(ii)]$\int_0^1 \int_0^1 \mathbb{E}|\mu(t,u)|^2 \,du \,dt < \infty$, where
$\int_0^1 \mathbb{E}|\mu(t,u)|^2 \,du < \infty$ for all $t\in[0,1]$.
\item[(iii)]$(\int_0^1 \mu(t,u) \,dW_u)_{t\in[0,1]}$ is a measurable modification.
\end{longlist}
Then, for $1< p < \infty$,
there exists a constant $c_{\fontsize{8.36pt}{10pt}\selectfont{(\ref{lemmadoubleIntBDG} )}}
\geq1 $ depending only on $p$
such that
\[
\biggl\llVert\biggl( \int_0^1 \biggl\llvert
\int_0^1 \mu(t,u) \,d W_u \biggr
\rrvert^2 \,dt \biggr) ^{1/2} \biggr\rrVert
_{L_p} \sim_{c_{\fontsize{6.6pt}{10pt}\selectfont{(\ref{lemmadoubleIntBDG} )}}} \biggl\llVert\biggl(
\int
_0^1 \int_0^1
\bigl| \mu(t,u) \bigr| ^2 \,du \,dt \biggr)^{1/2} \biggr\rrVert
_{L_p}.
\]
\end{lemma}

\begin{pf}
For the convenience of the reader, we sketch the proof. By a further
modification, we can assume that
$((\int_0^1 \mu(t,u) \,dW_u)(\omega))_{t\in[0,1]}\in L_2[0,1]$ for
all $\omega\in\Omega$
because of assumption (ii).
Assume that $(h_n)_{n=0}^\infty$ is the orthonormal basis of
Haar-functions in $L_2[0,1]$
and that $\mu_k(t,u)$ is the $k$th row of $\mu(t,u)$. Letting
\[
L_u^{n,k}:= \int_0^1
h_n(t) \mu_k(t,u) \,dt
\]
and using a stochastic Fubini argument we see that
\[
\biggl\llVert\biggl( \int_0^1 \biggl\llvert
\int_0^1 \mu(t,u) \,d W_u \biggr
\rrvert^2 \,dt \biggr) ^{1/2} \biggr\rrVert
_{L_p} = \Biggl\llVert\Biggl( \sum_{n=0}^\infty
\sum_{k=1}^d \biggl\llvert\int
_0^1 L_u^{n,k}
\,dW_u \biggr\rrvert^2 \Biggr) ^{1/2}
\Biggr\rrVert_{L_p}.
\]
Using the Burkholder--Davis--Gundy inequalities (\ref{eqnBDG}), we
obtain that
\begin{eqnarray*}
\biggl\llVert\biggl( \int_0^1 \biggl\llvert
\int_0^1 \mu(t,u) \,d W_u \biggr
\rrvert^2 \,dt \biggr) ^{1/2} \biggr\rrVert
_{L_p} &\sim_{c_{\fontsize{6.6pt}{10pt}\selectfont{(\ref{eqnBDG})}}}& \Biggl\llVert\Biggl( \sum
_{n=0}^\infty\sum_{k=1}^d
\int_0^1 \bigl|L_u^{n,k}\bigr|^2
\,d u \Biggr) ^{1/2} \Biggr\rrVert_{L_p}
\\
& = & \Biggl\llVert\Biggl( \sum_{k=1}^d
\int_0^1 \int_0^1
\bigl|\mu_k(t,u)\bigr|^2 \,dt \,d u \Biggr) ^{1/2}
\Biggr\rrVert_{L_p}
\\
& = & \biggl\llVert\biggl( \int_0^1 \int
_0^1 \bigl|\mu(t,u)\bigr|^2 \,dt \,d u \biggr)
^{1/2} \biggr\rrVert_{L_p}.
\end{eqnarray*}\upqed
\end{pf}

%
\begin{lemma}\label{lemmahypercontraction}
Let $1<p<\infty$, $g(Y_1)\in L_p$ and $0<t<1$,
and let $a=(a_1,\ldots,a_d)$ be a multi-index of differentiation.
Assume that $G$ is given by (\ref{eqndefinitionG}). Then
\[
\Bigl\llVert\sup_{0\le s \le t} \bigl|D_y^a
G(s,Y_s)\bigr| \Bigr\rrVert_{L_q}<\infty\qquad\mbox{for }
0<q< q(p,t):= \frac{p-1+t}{t}.
\]
\end{lemma}

\begin{pf*}{Sketch of the proof}
We use the notation (\ref{eqnconv1}) and (\ref{eqnconv2}) and
consider first the case
that $Y$ is the Brownian motion. A simple direct computation gives the
hyper-contraction property
\[
\bigl|D^a_x F(t,x)\bigr| \le C(q,t,a) \| f \|_{L_p(\mathbb{R}^d,\gamma_d)}
e^{|x|^2/(2tq)}
\]
for $0 < t < 1$ and $0<q<q(p,t)$. Moreover, the identity
\[
D^a_xF(s,x) = \mathbb{E}D^a_x
F (t,x+W_{t-s})
\]
for $0\le s \le t < 1$ directly implies that $(D^a_s F(s,W_s))_{s\in[0,t]}$
is an $L_q$-martingale. Therefore, we can exploit Doob's maximal inequality
for $1<q<q(p,t)$ to conclude
%
\begin{equation}
\label{eqnhypermaximalW} \mathbb{E}\sup_{0\le s \le t} \bigl|D^a_x
F(s,W_s)\bigr|^q <\infty\qquad\mbox{for all }0<q<
q(p,t).
\end{equation}
The case of the geometric Brownian motion can be deduced from the case
of the Brownian motion.
Using the notation (\ref{eqnconv1}) and (\ref{eqnconv2}) to
switch between the
Brownian motion and the geometric Brownian motion, we get for $0\le t <
1$ that
\[
D^a_y G(t,Y_t) = \Biggl[ \prod
_{k=1}^d \bigl(Y_t^k
\bigr)^{-a_k} \Biggr] \sum_{0\le b \le a}
\kappa_a^b D^b_xF(t,W_t),
\]
where $0\le b \le a$ is the coordinate-wise ordering and $\kappa_a^b$
are fixed coefficients.
Using (\ref{eqnhypermaximalW}), the integrability properties of the
geometric
Brownian motion and H\"older's inequality, we conclude that
\[
\label{eqnhypermaximalS} \mathbb{E}\sup_{0\le s \le t} \bigl|D^a_y
G(s,Y_s)\bigr|^q <\infty\qquad\mbox{for all }0<q<
q(p,t).
\]\upqed
\end{pf*}

The following estimates are known for more general processes than
the Brownian motion (see
\cite{Gobet-Munos2005} and
\cite{Geiss-Geiss-Gobet2012}, Remark 3). In our case,
they can be easily verified by using the martingale property
of the processes $(\nabla F(t,W_t))_{t \in[0,1)}$ and
$(D^2 F(t,W_t))_{t\in[0,1)}$.

%
\begin{lemma} \label{lemmaestimateDerFbyf}
Let $2\le p < \infty$. Assume that $f\dvtx \mathbb{R}^d\to\mathbb{R}$ is
measurable with
$f\in L_p(\mathbb{R}^d,\gamma_d)$ and that $F\dvtx [0,1]\times\mathbb
{R}^d \to\mathbb{R}$ is given
by $F(t,x):= \mathbb{E}f(x+W_{1-t})$.
Then there exists a constant $c_{\fontsize{8.36pt}{10pt}\selectfont{(\ref{lemmaestimateDerFbyf} )}} > 0$
depending only on $p$ such that, for all $0\le t<1$,
\begin{longlist}[(ii)]
\item[(i)]
$ \|\nabla F(t, W_t) \|_{L_{p}} \leq
c_{\fontsize{8.36pt}{10pt}\selectfont{(\ref{lemmaestimateDerFbyf} )}} (1-t)^{-1/2}
\| f(W_1) - F(t,W_t)\|_{L_{p}} $,
\item[(ii)]
$ \| D^2 F(t, W_t) \|_{L_{p}} \leq
c_{\fontsize{8.36pt}{10pt}\selectfont{(\ref{lemmaestimateDerFbyf} )}} (1-t)^{-1} \| f(W_1)-F(t,W_t) \|_{L_{p}} $.
\end{longlist}
\end{lemma}

Next, we state some Hardy type inequalities we have used in the paper.

%
\begin{proposition}\label{propgeneralinterpol}
Let $0 < \theta< 1$, $2 \leq q \leq\infty$ and let
$d^k\dvtx [0,1)\to[0,\infty)$, $k=0,1,2$, be measurable functions.
Assume that
\[
\frac{1}{\alpha} (1-t)^{k/2} \,d^k (t) \le
d^0(t) \le\alpha\biggl( \int_t^1
\bigl[d^1(s)\bigr]^2 \,ds \biggr) ^{1/2}
\qquad\mbox{for } t\in[0,1)
\]
and for $k=1,2$, and that
\[
d^1(t) \le A + \alpha\biggl( \int_0^t
\bigl[d^2(u)\bigr]^2 \,du \biggr) ^{1/2}
\qquad\mbox{for } t\in[0,1)
\]
for some $A\ge0$ and $\alpha>0$. Then
\begin{eqnarray*}
&& \bigl\|(1-t)^{-\theta/2} \,d^0(t) \bigr\|_{L_q ( [0,1), (dt)/(1-t) ) }
\\
&&\qquad  \sim_{c_{\fontsize{6.6pt}{10pt}\selectfont{(\ref{propgeneralinterpol})}}}\bigl\|
(1-t)^{(1-\theta)/2} \,d^1(t) \bigr\|
_{L_q ( [0,1), (dt)/(1-t) ) }
\end{eqnarray*}
and
\begin{eqnarray*}
& & \bigl\|(1-t)^{(2-\theta)/2} \,d^2(t) \bigr\|_{L_q ( [0,1), (dt)/(1-t) ) }
\\
&&\qquad \le c_{\fontsize{8.36pt}{10pt}\selectfont{(\ref{propgeneralinterpol})}} \bigl\|
(1-t)^{-\theta/2}
\,d^0(t) \bigr\|_{L_q ( [0,1), (dt)/(1-t) ) }
\\
&&\qquad \le c_{\fontsize{8.36pt}{10pt}\selectfont{(\ref{propgeneralinterpol})}}^2 \bigl[ A + \bigl\|
(1-t)^{(2-\theta)/2} \,d^2(t) \bigr\|
_{L_q ( [0,1), (dt)/(1-t) ) } \bigr],
\end{eqnarray*}
where $c_{\fontsize{8.36pt}{10pt}\selectfont{(\ref{propgeneralinterpol})}}\ge1$ depends at most on
$(\alpha,\theta,q)$.
If the functions $d^1$ and $d^2$ are nondecreasing, then the
inequalities are true for
$1\le q < 2$ as well.
\end{proposition}

From Proposition~\ref{propgeneralinterpol}, it follows that
%
\begin{equation}
\label{eqngeneralinterpolcor} A + \bigl\|
(1-t)^{(k-\theta)/2}
\,d^k(t) \bigr\|_{L_q}
\sim_{c_{\fontsize{6.6pt}{10pt}\selectfont{(\ref{eqngeneralinterpolcor})}}} A + \bigl\| (1-t)^{(l-\theta)/2}
\,d^l(t) \bigr\|_{L_q}
\end{equation}
for $L_q=L_q ( [0,1), \frac{dt}{1-t} ) $, $k,l=0,1,2$
and $c_{\fontsize{8.36pt}{10pt}\selectfont{(\ref{eqngeneralinterpolcor})}}:= [1 + c_{\fontsize{8.36pt}{10pt}\selectfont{(\ref{propgeneralinterpol})}}]^2$.
To prove Proposition~\ref{propgeneralinterpol}, we need:

%
\begin{lemma} \label{lemmaNonHardyIneq}
Let $0 < \theta< 1$, $2\le q \le\infty$ and let $\phi\dvtx  [0,1) \to
[0, \infty)$ be a measurable function.
Then there is a constant $c_{\fontsize{8.36pt}{10pt}\selectfont{(\ref{lemmaNonHardyIneq} )}}
> 0$, depending at most on $\theta$, such that
%
\begin{eqnarray}
&& \biggl\|(1-t)^{(1-\theta)/2} \biggl( \int
_0^t \phi(u)^2 \,du
\biggr)^{1/2} \biggr\|_{L_{q} ([0,1), (dt)/(1-t) )}
\nonumber\\[-8pt]\\[-8pt]
&&\qquad \leq c_{\fontsize{8.36pt}{10pt}\selectfont{(\ref{lemmaNonHardyIneq} )}} \bigl\|
(1-t)^{1-(\theta/2)} \phi(t)
\bigr\|_{L_{q} ([0,1), (dt)/(1-t))}.\nonumber
\end{eqnarray}
Moreover, if $\phi$ is nondecreasing, the inequality is true
for $1 \leq q < 2$ as well.
\end{lemma}

\begin{pf}
(a) For $2 \leq q \leq\infty$, we can use Hardy's inequality (see,
e.g., \cite{B-S}, Theorem 3.3.9):
for $- \infty< \lambda< 1$ and $1 \leq r < \infty$, and a
measurable $\psi\dvtx  (0, \infty) \to[0, \infty)$,
\[
\biggl( \int_0^{\infty} \biggl[t^{1-\lambda}
\int_t^{\infty} \psi(s) \frac{ds}{s}
\biggr]^r \frac{dt}{t} \biggr)^{1/r} \leq
\frac{1}{1-\lambda} \biggl( \int_0^{\infty}
\bigl[t^{1-\lambda} \psi(t) \bigr]^r \frac{dt}{t}
\biggr)^{1/r}
\]
and the same with the supremum norm if $r=\infty$.
With the notation $r:= \frac{q}{2}$, $g(t)=[\phi(t)]^2$, and $
\lambda=\theta$, we compute,
in the case $2\le q < \infty$,
\begin{eqnarray*}
&& \biggl\|(1-t)^{(1-\theta)/2} \biggl( \int
_0^t \phi(u)^2 \,du
\biggr)^{1/2} \biggr\|_{L_{q}
([0,1), (dt)/(1-t))}^2
\\
&&\qquad  =  \biggl(\int_0^1 \biggl[
(1-t)^{1-\theta} \int_0^t g(u) \,du
\biggr]^r \frac{dt}{1-t} \biggr)^{1/r}
\\
&&\qquad  =  \biggl(\int_0^{\infty} \biggl[
s^{1-\theta} \int_s^{\infty} h(v) \,dv
\biggr]^r \frac{ds}{s} \biggr)^{1/r},
\end{eqnarray*}
where $h(v) = g(1-v) \chi_{(0,1]} ( v )$.
Now we use Hardy's inequality for $\psi(v) = v h(v)$ and continue with
\begin{eqnarray*}
&& \biggl(\int_0^{\infty} \biggl[ s^{1-\theta}
\int_s^{\infty} \psi(v) \frac{dv}{v}
\biggr]^r \frac{ds}{s} \biggr)^{1/r}
\\
&&\qquad  \leq
\frac{1}{1-\theta} \biggl(\int_0^{\infty} \bigl[
s^{1-\theta} \psi(s) \bigr]^r \frac{ds}{s}
\biggr)^{1/r}
\\
&&\qquad  =  \frac{1}{1-\theta} \biggl(\int_0^{\infty}
\bigl[ s^{2-\theta} h(s) \bigr]^r \frac{ds}{s}
\biggr)^{1/r}
\\
&&\qquad  =  \frac{1}{1-\theta} \bigl\|(1-t)^{1-(\theta/2)} \phi(t)
\bigr\|_{L_{q} ([0,1), (dt)/(1-t))}^2
\end{eqnarray*}
and the proof is complete for
$2 \leq q < \infty$. The case $q= \infty$ is analogous.

(b) For $1 \leq q < 2,$ we use a different argument.
First, we define $r:= \frac{2}{q}$ so that $1 < r \leq2$.
For $0<T<1,$ we compute
\begin{eqnarray*}
&& \int_0^1 (1-t)^{(1-\theta)/r} \biggl( \int
_0^t \chi_{ [
T,1 )} ( u ) \,du
\biggr)^{1/r} \frac{dt}{1-t}
\\
&&\qquad  = \int_0^1
(1-t)^{(1-\theta)/r} (t-T)_+^{1/r} \frac{dt}{1-t}
\\
&&\qquad \leq (1-T)^{1/r} \int_T^1
(1-t)^{((1-\theta)/r) -1} \,dt
\\
&&\qquad  =  c \int_T^1 (1-t)^{(2-\theta)/r}
\chi_{ [ T,1
)} ( t ) \frac{dt}{1-t}
\end{eqnarray*}
with $c:= \frac{2-\theta}{1-\theta}$.
This proves the desired inequality for $\psi^{(T)} (t):= \chi_{
[ T,1 )} ( t )$.
Next, we define $\psi:=\phi^q$ so that $\psi^r = \phi^2$.
By assumption, $\phi$ is nondecreasing, and so is~$\psi$, too.
Now,\vspace*{1pt} we can approximate $\psi$ from below by a sum of functions like
$\psi^{(T)}$:
for each integer $n \geq1$, we find $\alpha_k^n \ge0$, $k = 0,\ldots, 2^n-1$ and
$0=t^n_0 < t^n_1 < \cdots< t^n_{2^n-1} < t^n_{2^n}=1$ such that
\[
\psi_n (t):= \sum_{k=0}^{2^n-1}
\alpha_k^n \psi^{(t^n_k)} (t) \to\psi(t)
\]
for almost all $t \in[0,1)$ and $\psi_{n-1} \leq\psi_n $ for
all $n\geq2$.
Then, since $r\geq1$,
\begin{eqnarray*}
& & \int_0^1 (1-t)^{(1-\theta)/r} \biggl(
\int_0^t \psi_n(u)^r
\,du \biggr)^{1/r} \frac{dt}{1-t}
\\
&&\qquad \leq \int_0^1 (1-t)^{(1-\theta)/r} \sum
_{k=0}^{2^n-1} \alpha_k^n
\biggl( \int_0^t \psi^{(t_k^n)} (u) \,du
\biggr)^{1/r} \frac{dt}{1-t}
\\
&&\qquad \leq \sum_{k=0}^{2^n-1}
\alpha_k^n c \int_0^1
(1-t)^{(2-\theta)/r} \psi^{(t_k^n)} (t) \frac{dt}{1-t}
\\
&&\qquad  =  c \int_0^1 (1-t)^{(2-\theta)/r}
\psi_n (t) \frac{dt}{1-t}
\end{eqnarray*}
and the claim follows by monotone convergence.
\end{pf}

\begin{pf*}{Proof of Proposition~\ref{propgeneralinterpol}}
(a) Our assumptions imply for all $1\le q \le\infty$ that
\[
\bigl\|(1-t)^{(1-\theta)/2} \,d^1(t) \bigr\|_{L_q ( [0,1), (dt)/(1-t) ) } \le\alpha\bigl\|
(1-t)^{-\theta/2} \,d^0(t) \bigr\|
_{L_q ( [0,1), (dt)/(1-t) ) }
\]
and
\[
\bigl\|(1-t)^{(2-\theta)/2} \,d^2(t) \bigr\|
_{L_q ( [0,1), (dt)/(1-t) ) } \le\alpha\bigl\|
(1-t)^{-\theta/2} \,d^0(t) \bigr\|
_{L_q ( [0,1), (dt)/(1-t) ) }.
\]

(b) Next, we observe that
\begin{eqnarray*}
& & \bigl\|(1-t)^{-\theta/2} \,d^0(t) \bigr\|_{L_q ( [0,1), (dt)/(1-t) ) }
\\
&&\qquad \le \alpha\biggl\|(1-t)^{(1-\theta)/2}
\biggl(
\frac{1}{1-t}\int_t^1 \bigl[d^1(s)
\bigr]^2 \,ds \biggr) ^{1/2} \biggr\|
_{L_q ( [0,1), (dt)/(1-t) ) }
\\
&&\qquad \le \alpha\theta^{-\max\{ 1/2,1/q \}} \bigl\|(1-t)^{(1-\theta)/2}
\,d^1(t) \bigr\|_{L_q ( [0,1), (dt)/(1-t) ) },
\end{eqnarray*}
where we used \cite{Gei19}, formula (14) (the condition that $\psi$
in \cite{Gei19} is continuous
in the case $1\le q < 2$ is not necessary).

(c) To prove the remaining inequality, we continue from (b)
with Lemma~\ref{lemmaNonHardyIneq} to
\begin{eqnarray*}
& & \bigl\|(1-t)^{(1-\theta)/2} \,d^1(t) \bigr\|_{L_q ( [0,1), (dt)/(1-t) ) }
\\
&&\qquad \le \biggl\|(1-t)^{(1-\theta)/2} \biggl[ A
+ \alpha
\biggl( \int_0^t d^2(u) \,du \biggr)
^{1/2} \biggr] \biggr\|_{L_q ( [0,1), (dt)/(1-t)) }
\\
&&\qquad \le A \bigl\|(1-t)^{(1-\theta)/2} \bigr\|_{L_q ( [0,1), (dt)/(1-t) ) }
\\
&&\quad\qquad{}+ \alpha c_{\fontsize{8.36pt}{10pt}\selectfont{(\ref{lemmaNonHardyIneq} )}} \bigl\|
(1-t)^{1-(\theta/2)} \,d^2(t) \bigr\|
_{L_q ( [0,1), (dt)/(1-t) ) }.
\end{eqnarray*}\upqed
\end{pf*}
\end{appendix}



%
%

\printaddresses

\end{document}